\documentclass[12pt]{amsart}
\usepackage{graphicx,amssymb}
\usepackage[all]{xy}
\usepackage{times}

\newtheorem{theorem}{Theorem}[section]
\newtheorem{proposition}[theorem]{Proposition}
\newtheorem{lemma}[theorem]{Lemma}

\theoremstyle{definition}
\newtheorem{definition}[theorem]{Definition}
\newtheorem{example}[theorem]{Example}

\theoremstyle{remark}
\newtheorem{remark}[theorem]{Remark}

\numberwithin{equation}{section}

\begin{document}

\title[Higher-order Alexander invariants for homologically fibered knots]
{Factorization formulas and computations of higher-order Alexander 
invariants for homologically fibered knots} 

\date{\today}


\author{Hiroshi Goda}
\address{Department of Mathematics,
Tokyo University of Agriculture and Technology,
2-24-16 Naka-cho, Koganei,
Tokyo 184-8588, Japan}
\email{goda@cc.tuat.ac.jp}
\author{Takuya Sakasai}
\address{Department of Mathematics,
Tokyo Institute of Technology, 
2-12-1 Oh-okayama, Meguro-ku, 
Tokyo 152-8551, Japan}
\email{sakasai@math.titech.ac.jp}

\thanks{
The authors are partially supported
by Grant-in-Aid for Scientific Research,
(No.~21540071 and No.~21740044),
Ministry of Education, Science, 
Sports and Technology, Japan. 
}


\begin{abstract}
Homologically fibered knots are knots whose exteriors 
satisfy the same homological conditions as fibered knots. 
In our previous paper, we observed that 
for such a knot, higher-order Alexander invariants 
defined by Cochran, Harvey and Friedl are generally 
factorized into the part of the Magnus matrix and 
that of a certain Reidemeister torsion, both of which are 
known as invariants of homology cylinders over a surface. 
In this paper, 
we study more details of the invariants 
and give their concrete calculations 
after restricting to the case of 
the invariants associated with 
metabelian quotients of knot groups. 
We provide explicit computational results 
of the invariants for all the 12 crossings 
non-fibered homologically fibered knots.
\end{abstract}

\maketitle


\section{Introduction}\label{sec:intro}

Let $K$ be a knot in a 3-sphere $S^3$. In our previous paper \cite{gs08}, 
we introduced a class of knots called {\it $($rationally$)$ 
homologically fibered knots} and 
studied their fundamental properties by using their Alexander invariants. 
A (rationally) homologically fibered knot $K$ is 
by definition a knot satisfying the property that the sutured manifold $M_R$ 
obtained from the exterior $E(K)$ of 
$K$ by cutting along a minimal 
genus Seifert surface $R$ is a (rational) homology product whose boundary is 
the union of two copies of $R$. 

For a rationally homologically fibered knot $K$ with a minimal genus 
Seifert surface $R$ of genus $g$, 
let $i_+, i_-:R \to \partial M_R$ denote the natural 
identifications of $R$ with 
the two sides of the boundary of $M_R$. 
We fix a basis of $H_1 (R;\mathbb{Q})$ giving rise to 
an isomorphism $H_1 (R;\mathbb{Q}) \cong \mathbb{Q}^{2g}$. 
Then, by 
using the invertibility (over $\mathbb{Q}$) of 
the Seifert matrix $S$, we can rewrite the definition 
$\Delta_K (t)=\det(S-tS^T)$ of the Alexander polynomial of $K$  
and obtain a factorization 
\begin{equation}\label{eq1}
\Delta_K(t) = \det (S) \det (I_{2g}-t \,\sigma (M_R))
\end{equation}
\noindent
of $\Delta_K (t)$. Here $\sigma(M_R):=S^{-1}S^T$ coincides with 
the representation matrix of the composite of isomorphisms
\[\mathbb{Q}^{2g} \cong H_1 (R;\mathbb{Q}) \xrightarrow[i_-]{\cong} 
H_1 (M_R;\mathbb{Q}) 
\xrightarrow[i_+^{-1}]{\cong} H_1 (R;\mathbb{Q}) \cong \mathbb{Q}^{2g}.\]
The matrix $\sigma(M_R)$ can be interpreted as a monodromy of $M_R$ 
from a view point of the rational homology. 
Regarding the formula $(\ref{eq1})$ as a basic case, 
in \cite{gs08} 
we gave its generalization 
under the framework of {\it higher-order Alexander invariants} due to 
Cochran \cite{coc}, Harvey \cite{har2} and Friedl \cite{fri}. 
In this procedure, the Seifert matrix $S$, the monodromy $\sigma (M_R)$ 
and $\Delta_K (t)$ are generalized to a certain Reidemeister torsion 
$\tau_\rho^+ (M_R)$, the Magnus matrix $r_\rho (M_R)$ and 
a higher-order (non-commutative) Alexander invariant $\tau_\rho (E(K))$ 
associated with a representation $\rho$ of 
the fundamental group of $M_R$. 
Then the generalized formula is given by 
\begin{equation}\label{eq_factor}
\tau_\rho (E(K)) = \frac{\tau_\rho^+ (M_R) \cdot 
(I_{2g} -\rho (\mu) r_\rho (M_R))}{1-\rho(\mu)}, 
\end{equation}
where $\mu \in \pi_1 (E(K))$ represents the meridian of $K$. 
To compare $(\ref{eq_factor})$ with $(\ref{eq1})$, 
recall Milnor's formula \cite{milnor} that 
$\displaystyle\frac{\Delta_K (t)}{1-t}$ 
represents a Reidemeister torsion associated with the abelianization 
homomorphism $\rho_1:\pi_1 (E(K)) \to \langle t \rangle$. 
For details of the formula, see Theorem \ref{thm:factorization}. 

The purpose of this paper is to investigate the factorization 
formula (\ref{eq_factor}) 
with explicit computational examples. 
In the theory of higher-order Alexander invariants, 
an 
important problem has been to find methods 
for computing the invariants and extract 
topological information from them. 
This problem arises from 
the difficulty in non-commutative rings involved in the definition. 
We now intend to understand the higher-order invariant $\tau_\rho (E(K))$ 
by looking at each of the constituents of the formula (\ref{eq_factor}).
More specifically, in the latter half of this paper, 
we focus on the invariants associated with 
metabelian quotients of knot groups of homologically fibered knots.  
In this situation, although $\tau_\rho (E(K))$ itself belongs to 
a non-commutative ring setting, both of $\tau_\rho^+ (M_R)$ and $r_\rho (M_R)$ 
can be computed in a realm of commutative rings. 
A sample calculation with details is given in Section \ref{sec:sample} 
and more examples are exhibited in Section \ref{sec:HFK12}, where 
we use $\tau_\rho^+ (M_R)$ to detect the non-fiberedness of all 
the 12 crossings non-fibered homologically fibered knots. 
We remark that in the situation of Sections \ref{sec:sample} 
and \ref{sec:HFK12}, the torsion 
$\tau_\rho^+ (M_R)$ may be regarded as 
a special case of a decategorification of 
the sutured Floer homology as shown by Friedl-Juh\'asz-Rasmussen \cite{fjr}. 
In Section \ref{sec:magnus}, 
we study the Magnus matrix $r_\rho (M_R)$ and 
see that $r_\rho (M_R)$ is unchanged under 
{\it concordances of Seifert surfaces} introduced by Myers \cite{myers}. 
Using his result, we mention how to obtain more examples of explicit 
computations. 

The authors would like to thank Professor Ko Honda for 
helpful discussions and Professor Robert Myers 
for informing the authors about 
his paper \cite{myers}. They also thank 
the anonymous referee for his-or-her 
helpful comments to improve 
the previous version of this paper.


\section{Homologically fibered knots and homology cylinders}\label{sec:HFK}

First, we recall the definition of 
sutured manifolds given by Gabai \cite{gabai1}. 
We here use a special case of them. 

A {\it sutured manifold\/} $(M,\gamma)$ is a compact 
oriented 3-manifold $M$ together with a subset $\gamma \subset \partial M$ 
which is a union of finitely many mutually disjoint annuli. 
For each component of $\gamma$, 
an oriented core circle called a {\it suture\/} is fixed, 
and we denote the set of sutures by $s(\gamma)$. 
Every component of $R(\gamma)=\partial M-{\rm Int\,}\gamma$
is oriented so that the orientations on $R(\gamma)$ are coherent 
with respect to $s(\gamma)$, that is, the orientation of each component 
of $\partial R(\gamma)$ induced from that of $R(\gamma)$ 
is parallel to the orientation of the corresponding component 
of $s(\gamma)$. 
We denote by $R_{+}(\gamma)$ (resp. $R_{-}(\gamma)$) 
the union of those 
components of $R(\gamma)$ whose normal vectors point out of 
(resp. into) $M$. 

\begin{example}
For a knot $K$ in $S^3$ and a Seifert surface $\bar{R}$ of $K$, 
we set $R:=\bar{R} \cap E(K)$, 
called also a Seifert surface, 
where $E(K)= \overline{S^3-N(K)}$ is the complement 
of a regular neighborhood $N(K)$ of $K$. 
Then $(M_R, \gamma):=(\overline{E(K)-N(R)}, 
\overline{\partial E(K)-N(\partial R)})$ 
defines a sutured manifold. 
We call it the {\it complementary sutured manifold} for $R$. 
In this paper, we simply call it the sutured manifold for $R$. 
\end{example}

\begin{figure}[htbp]
\begin{center}
\includegraphics[width=0.6\textwidth]{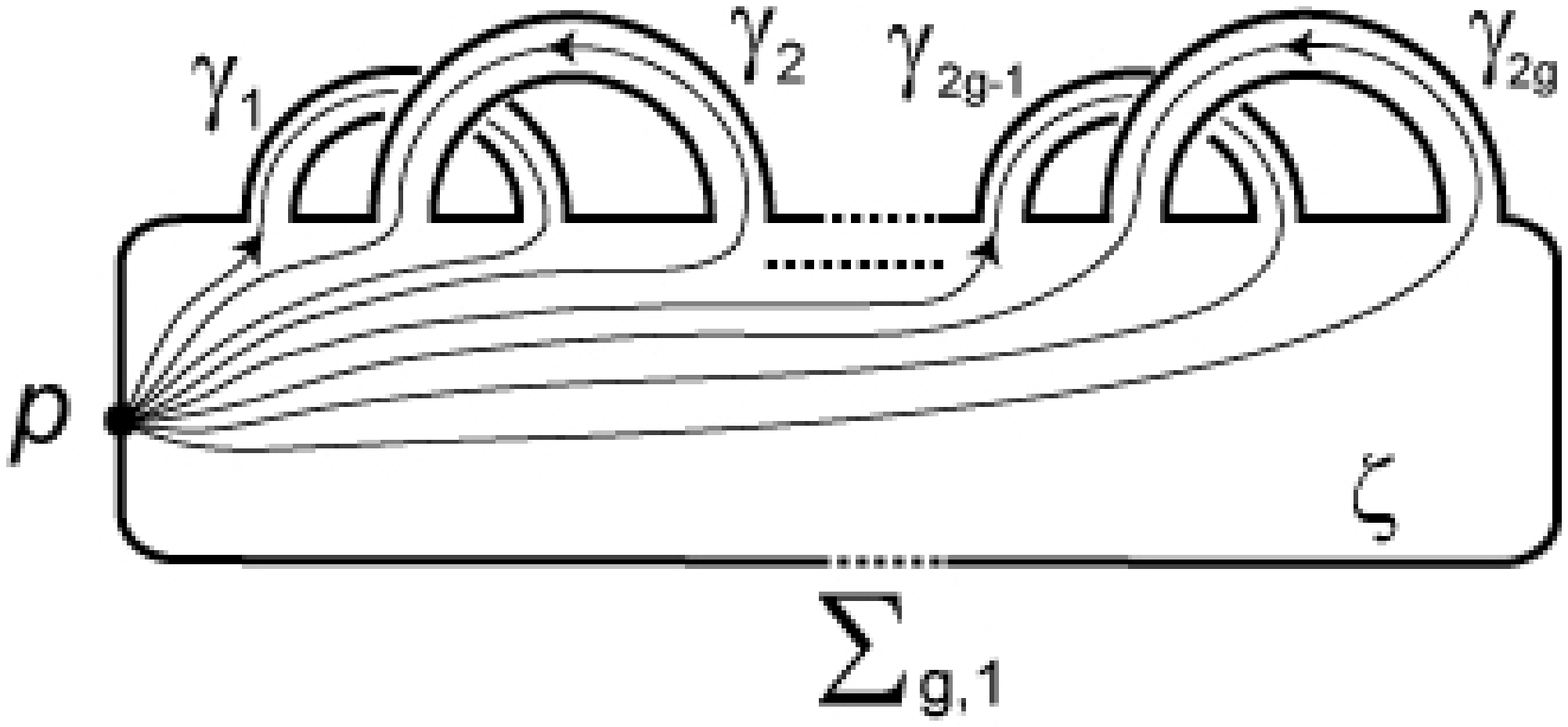}
\end{center}
\caption{Cell decomposition of $\Sigma_{g,1}$}
\label{fig:spine1}
\end{figure}

\begin{definition}[\cite{gs08}]\label{def:HFknot}
A knot $K$ in $S^3$ is called 
a {\it rationally homologically fibered knot} of genus $g$ if 
it has the following properties which are equivalent to each other: 
\begin{itemize}
\item[(a)] The degree of the Alexander polynomial $\Delta_K (t)$ of 
$K$ is equal to twice the genus $g=g(K)$ of $K$; 

\item[(b)] For any minimal genus Seifert surface $R$ of $K$, its 
Seifert matrix $S$ is invertible over $\mathbb{Q}$; and 

\item[(c)] The sutured manifold $(M_R,\gamma)$ 
for any minimal genus Seifert surface $R$ is a rational 
homology product over $R$. 
\end{itemize}
Moreover, when $\Delta_K (t)$ is monic 
(correspondingly, $S$ is invertible over $\mathbb{Z}$ 
and $M_R$ is a homology product), we say 
$K$ is a {\it homologically fibered knot}.
\end{definition}
\begin{remark}
Aside from the name, 
the equivalence of the conditions (a), (b), (c) 
in Definition \ref{def:HFknot} 
was mentioned in Crowell-Trotter \cite{ct}. 
\end{remark}
Next we recall the definition of {\it homology cylinders}, which 
can be regarded as a generalization of mapping classes of surfaces. 
We refer to Goussarov \cite{gou}, Habiro \cite{habiro}, 
Garoufalidis-Levine \cite{gl} and Levine \cite{levin} for 
their origin. Strictly speaking, the definition below is 
closer to that in \cite{gl} and \cite{levin}. 
Let $\Sigma_{g,1}$ be a compact connected 
oriented surface of genus $g\ge 0$ with a connected boundary. 
We fix a cell decomposition of $\Sigma_{g,1}$ 
consisting of one vertex $p$, 
edges $\gamma_1, \gamma_2, \ldots, \gamma_{2g}, \zeta$ and 
one face as in Figure \ref{fig:spine1}. 

\begin{definition}
A {\it homology cylinder\/}  $(M,i_{+},i_{-})$ {\it over} $\Sigma_{g,1}$ 
consists of a compact oriented 3-manifold $M$ 
with two embeddings 
$i_{+}, i_{-}: \Sigma_{g,1} \hookrightarrow \partial M$, 
called {\it markings\/}, 
such that:
\begin{enumerate}
\renewcommand{\labelenumi}{(\roman{enumi})}
\item
$i_{+}$ is orientation-preserving and $i_{-}$ is orientation-reversing; 
\item 
$\partial M=i_{+}(\Sigma_{g,1})\cup i_{-}(\Sigma_{g,1})$ and 
$i_{+}(\Sigma_{g,1})\cap i_{-}(\Sigma_{g,1})
=i_{+}(\partial\Sigma_{g,1})=i_{-}(\partial\Sigma_{g,1})$;
\item
$i_{+}|_{\partial \Sigma_{g,1}}=i_{-}|_{\partial \Sigma_{g,1}}$; and 
\item
$i_{+},i_{-} : H_{*}(\Sigma_{g,1};\mathbb Z)\to H_{*}(M;\mathbb Z)$ 
are isomorphisms. 
\end{enumerate}
Similarly, the definition of a {\it rational homology cylinder\/} is 
obtained by replacing (iv) with the condition that (iv') 
$i_{+},i_{-} : H_{*}(\Sigma_{g,1};\mathbb Q)\to H_{*}(M;\mathbb Q)$ 
are isomorphisms. 
\end{definition}
Two homology cylinders $(M,i_+,i_-)$ and $(N,j_+,j_-)$ over $\Sigma_{g,1}$ 
are said to be {\it isomorphic} if there exists 
an orientation-preserving diffeomorphism $f:M \xrightarrow{\cong} N$ 
satisfying $j_+ = f \circ i_+$ and $j_- = f \circ i_-$. 
We denote by $\mathcal{C}_{g,1}$ the set of all isomorphism classes 
of homology cylinders over $\Sigma_{g,1}$. 
By using markings, we can endow 
$\mathcal{C}_{g,1}$ with a monoid structure whose product is given by 
\[(M,i_+,i_-) \cdot (N,j_+,j_-)
:=(M \cup_{i_- \circ (j_+)^{-1}} N, i_+,j_-)\]
for $(M,i_+,i_-)$, $(N,j_+,j_-) \in \mathcal{C}_{g,1}$. 
The unit of this monoid is given by 
\[(M,i_+,i_-)=(\Sigma_{g,1} \times [0,1], \mathrm{id} 
\times 1, \mathrm{id} \times 0),\]
where collars of $i_+ (\Sigma_{g,1})$ and $i_- (\Sigma_{g,1})$ 
are stretched half-way along $(\partial \Sigma_{g,1}) \times [0,1]$. 
The monoid $\mathcal{C}_{g,1}^\mathbb{Q}$ 
of all isomorphism classes of rational homology cylinders 
over $\Sigma_{g,1}$ is defined similarly. 
For each diffeomorphism $\varphi$ of 
$\Sigma_{g,1}$ which fixes $\partial \Sigma_{g,1}$ pointwise, 
we can construct a homology cylinder as a {\it mapping cylinder} 
\[(\Sigma_{g,1} \times [0,1], \mathrm{id} \times 1, 
\varphi \times 0)\]
of $\varphi$. 

Constructing a homology cylinder from a given homologically fibered knot has 
an ambiguity arising from taking a minimal genus Seifert surface and fixing 
a pair of markings. 

\begin{proposition}\label{prop:unique}
Let $R_1$ and $R_2$ be $($maybe parallel\/$)$ 
minimal genus Seifert surfaces of a homologically 
fibered knot of genus $g$ 
and let $M_{R_1}$ and $M_{R_2}$ be their sutured 
manifolds. For any markings $i_\pm$ and $j_\pm$ of 
$\partial M_{R_1}$ and $\partial M_{R_2}$, 
there exists another 
homology cylinder $N \in \mathcal{C}_{g,1}$ such that 
\[(M_{R_1},i_+,i_-) \cdot N = N \cdot (M_{R_2},j_+,j_-)\]
holds as elements of $\mathcal{C}_{g,1}$. 
\end{proposition}
\begin{proof}
First we assume that $R_1$ and $R_2$ are disjoint in $E(K)$. 
Cut $E(K)$ along $R_1$ and $R_2$. Then we obtain two submanifolds 
$N$ and $N'$ of $E(K)$, where $N$ (resp. $N'$)
may be regarded as a surface cobordism 
between $i_+ (\Sigma_{g,1})$ and $j_- (\Sigma_{g,1})$ 
(resp. $j_+ (\Sigma_{g,1})$ and $i_- (\Sigma_{g,1})$). 
We can easily check that 
$(N,i_+,j_-)$ and $(N',j_+,i_-)$ are homology cylinders 
over $\Sigma_{g,1}$. Then the equality 
$M_{R_{1}} \cup_{R_1} N=N \cup_{R_2} M_{R_2}$ holds and it 
shows our claim in this case. 

For the general case, we can use a theorem of Scharlemann-Thompson \cite{st} . 
It says that there exists a  sequence of minimal genus Seifert surfaces 
$R_1=S_1\rightarrow S_2 \rightarrow \cdots \rightarrow S_n=R_2$ 
such that $S_i$ and $S_{i+1}$ are disjoint in $E(K)$
for $i=1,2,\ldots , n-1$. 
Using the above argument repeatedly, 
we have the conclusion.  
\end{proof}
\noindent
This proposition can be seen as a generalization of the fact that 
a fibered knot determines an element of the mapping class group of 
a surface uniquely up to conjugation. 
\begin{remark}
Differently from fibered knots, a homologically fibered knot does not 
necessarily have a unique minimal genus Seifert surface. 
Indeed, it was shown by Eisner \cite{eisner} that 
the connected sum of two non-fibered knots has infinitely many 
non-isotopic minimal genus Seifert surfaces. Hence the connected sum of 
two non-fibered homologically fibered knots, which is again a homologically 
fibered knot, gives such an example. 
The authors do not know whether there exists a homologically fibered knot 
which has minimal genus Seifert surfaces whose complements are 
not homeomorphic. 
\end{remark}


\section{Higher-order Alexander invariants}\label{sec:alexander}

From the factorization $(\ref{eq1})$, we see 
that if a rationally homologically fibered 
knot has a non-trivial $\det (S)$-part, that is $|\det (S)| \neq 1$, then 
this knot is not fibered. However, this argument is useless 
for homologically fibered knots, since $|\det (S)|= 1$. 
In this section, we give a generalization of the factorization $(\ref{eq1})$ 
by using the framework of {\it higher-order Alexander invariants} originally 
due to Cochran \cite{coc} and Harvey \cite{har2} together with 
their interpretations as Reidemeister torsions given by Friedl \cite{fri}. 
We will see later that this generalized factorization works 
well for homologically fibered knots. 

We begin by summarizing our notation. 
For a matrix $A$ with entries in a group ring 
$\mathbb{Z} G$ (or its quotient field) 
for a group $G$, we denote by $\overline{A}$ 
the matrix obtained from $A$ by applying the involution induced 
from $(x \mapsto x^{-1},\ x \in G)$ to each entry. 
For a module $M$, we write $M^n$ for the module of 
column vectors with $n$ entries. 
For a finite cell complex $X$, we denote by $\widetilde{X}$ 
its universal covering. We take a base point $p$ of $X$ and 
a lift $\widetilde{p}$ of $p$ as a base point of $\widetilde{X}$. 
$\pi:=\pi_1 (X,p)$ acts on $\widetilde{X}$ from 
the {\it right} through its deck transformation group, so that 
the lift of a loop $l \in \pi$ starting from $\widetilde{p}$ 
reaches $\widetilde{p} \, l^{-1}$. 
Then the cellular chain complex $C_{\ast} (\widetilde{X})$ of 
$\widetilde{X}$ 
becomes a right $\mathbb{Z} \pi$-module. 
For each left $\mathbb{Z}\pi$-algebra $\mathcal{R}$, 
the twisted chain complex $C_{\ast} (X;\mathcal{R})$ is given 
by the tensor product of 
the right $\mathbb{Z}\pi$-module $C_{\ast} (\widetilde{X})$ and 
the left $\mathbb{Z}\pi$-module $\mathcal{R}$, 
so that $C_{\ast} (X;\mathcal{R})$ and 
$H_{\ast} (X;\mathcal{R})$ are right $\mathcal{R}$-modules. 

In the definition of higher-order Alexander invariants, PTFA groups play 
important roles, where a group $\Gamma$ is said 
to be {\it poly-torsion-free abelian $($PTFA$)$} 
if it has a sequence 
\[\Gamma=\Gamma_0 \triangleright \Gamma_1 \triangleright 
\cdots \triangleright \Gamma_n = \{1\}\]
whose successive quotients $\Gamma_i/\Gamma_{i+1}$ $(i \ge 0)$ 
are all torsion-free abelian. An advantage of using PTFA groups is that 
the group ring $\mathbb{Z} \Gamma$ (or $\mathbb{Q} \Gamma$) of $\Gamma$ 
is known to be an {\it Ore domain} so that it can be embed into 
the field (skew field in general) 
\[\mathcal{K}_\Gamma:= \mathbb{Z}\Gamma (\mathbb{Z}\Gamma - \{0\})^{-1}
=\mathbb{Q}\Gamma (\mathbb{Q}\Gamma - \{0\})^{-1}\]
called the {\it right field of fractions}. 
We refer to \cite{coc} and \cite{pa} for generalities of 
PTFA groups and localizations of their group rings. 
A typical example of PTFA groups is $\mathbb{Z}^n$, 
where $\mathcal K_{\mathbb{Z}^n}$ is isomorphic to 
the field of rational functions with $n$ variables.

For a rationally homologically fibered knot $K$, we take 
a {\it non-trivial} homomorphism 
$\rho:G(K) \to \Gamma$ to a PTFA group $\Gamma$, 
where $G(K)$ denotes the knot group $\pi_1 (E(K))$. 
We can regard $\mathcal{K}_\Gamma$ as a local coefficient system 
on $E(K)$ through $\rho$. 
Using arguments in Cochran-Orr-Teichner \cite[Section 2]{cot} and 
Cochran \cite[Section 3]{coc}, we have: 
\begin{lemma}\label{lem:vanish}
For any non-trivial homomorphism $\rho:G(K) \to \Gamma$ to 
a PTFA group $\Gamma$, 
we have $H_\ast (E(K); \mathcal{K}_\Gamma)=0$.
%
\end{lemma}
\noindent
By this lemma, 
we can define the Reidemeister torsion 
\[\tau_{\rho} (E(K)):=
\tau(C_\ast (E(K); \mathcal{K}_{\Gamma})) 
\in K_1(\mathcal{K}_\Gamma)/\pm\rho(G(K))\] 
for the acyclic complex $C_\ast (E(K); \mathcal{K}_{\Gamma})$. 
We refer to Milnor \cite{milnor} 
for generalities of torsions. 
By higher-order Alexander invariants for $K$, 
we here mean this torsion $\tau_{\rho} (E(K))$. 

We now describe a factorization of $\tau_{\rho} (E(K))$ 
generalizing (\ref{eq1}). 
For that we use two kinds of invariants for rational 
homology cylinders from \cite{sakasai08} and \cite{gs08}. 
Let $(M_R,i_+,i_-) \in \mathcal{C}_{g,1}^\mathbb{Q}$ 
be the rational homology cylinder 
obtained as the sutured manifold for 
a minimal genus Seifert surface $R$ of $K$. 
We use the same notation $\rho: \pi_1 (M_R) \to \Gamma$ for the 
composition $\pi_1 (M_R) \to G(K) \xrightarrow{\rho} \Gamma$. 
By applying Cochran-Orr-Teichner \cite[Proposition 2.10]{cot}, 
we have: 
\begin{lemma}\label{lem:vanish2}
$i_+, i_-: H_\ast (\Sigma_{g,1},p ;i_{\pm}^\ast \mathcal{K}_\Gamma) 
\to H_\ast (M_R,p ;\mathcal{K}_\Gamma)$ are isomorphisms 
as right $\mathcal{K}_\Gamma$-vector spaces. Equivalently, 
$H_\ast (M_R,i_\pm(\Sigma_{g,1}); \mathcal{K}_\Gamma)= 0$. 
\end{lemma}
\noindent
This lemma provides the following two kinds of invariants for $M_R$.  

\bigskip
\noindent
{\bf The Magnus matrix} \ 
Let $X \subset \Sigma_{g,1}$ be the union of $2g$ loops 
$\gamma_1, \ldots , \gamma_{2g}$ 
(see Figure \ref{fig:spine1}). 
$X$ is a deformation retract of $\Sigma_{g,1}$ relative to $p$. 
Therefore, for $\pm \in \{+,-\}$, we have
\[H_1 (\Sigma_{g,1},p;i_{\pm}^\ast \mathcal{K}_\Gamma) \cong 
H_1 (X,p;i_{\pm}^\ast \mathcal{K}_\Gamma) = 
C_1 (\widetilde{X}) \otimes_{\pi_1 (\Sigma_{g,1})} 
i_{\pm}^\ast \mathcal{K}_\Gamma \cong \mathcal{K}_\Gamma^{2g}\]
with a basis
\[\{ \widetilde{\gamma}_1 \otimes 1, \ldots , 
\widetilde{\gamma}_{2g} \otimes 1\} 
\subset C_1 (\widetilde{X}) \otimes_{\pi_1 (\Sigma_{g,1})} 
i_{\pm}^\ast \mathcal{K}_\Gamma\] 
as a right $\mathcal{K}_\Gamma$-vector space. 
Here we fix a lift $\widetilde{p}$ of $p$ as a 
base point of $\widetilde{X}$, and 
denote by $\widetilde{\gamma}_i$ the lift of 
the oriented loop $\gamma_i$ 
starting from $\widetilde{p}$. 

\begin{definition}\label{def:Mag2}
For $M_R=(M_R,i_+,i_-) \in \mathcal{C}_{g,1}^\mathbb{Q}$, 
the {\it Magnus matrix} 
\[r_\rho (M_R) \in GL(2g,\mathcal{K}_\Gamma)\]
of $M_R$ is defined as the representation matrix of 
the right $\mathcal{K}_\Gamma$-isomorphism
\[\mathcal{K}_\Gamma^{2g} \cong 
H_1 (\Sigma_{g,1},p;i_-^\ast \mathcal{K}_\Gamma) \xrightarrow[i_-]{\cong} 
H_1 (M_R,p;\mathcal{K}_\Gamma) \xrightarrow[i_+^{-1}]{\cong} 
H_1 (\Sigma_{g,1},p;i_+^\ast \mathcal{K}_\Gamma) 
\cong \mathcal{K}_\Gamma^{2g},\]
where the first and the last isomorphisms use 
the bases mentioned above.
\end{definition}
\noindent
The matrix $r_\rho(M_R)$ can be interpreted as a monodromy of $M_R$ 
from a view point of the twisted homology with coefficients 
in $\mathcal{K}_\Gamma$. 

\bigskip

\noindent
{\bf $\Gamma$-torsion} \ 
Since the relative complex 
$C_\ast (M_R,i_+(\Sigma_{g,1});\mathcal{K}_\Gamma)$ 
obtained from any cell decomposition of $(M_R, i_+(\Sigma_{g,1}))$ 
is acyclic by Lemma \ref{lem:vanish2}, 
we can define the following:
\begin{definition}
For $M_R=(M_R,i_+,i_-) \in \mathcal{C}_{g,1}^\mathbb{Q}$, 
the $\Gamma$-{\it torsion} $\tau_{\rho}^+ (M_R)$ 
of $M_R$ is defined by 
\[\tau_{\rho}^+ (M_R):=
\tau(C_\ast (M_R,i_+(\Sigma_{g,1});\mathcal{K}_\Gamma))
\in K_1 (\mathcal{K}_\Gamma)/\pm \rho (\pi_1 (M_R)).\]
\end{definition}

\bigskip
A method for computing 
$r_\rho(M_R)$ and $\tau_{\rho}^+ (M_R)$ is given in 
\cite[Section 4]{gs08}, which is based on 
Kirk-Livingston-Wang's method \cite{klw} for 
invariants of string links, and we now recall it briefly. 
An {\it admissible presentation} of $\pi_1 (M_R)$ is defined to be 
the one of the form 
\begin{align}\label{admissible}
\langle i_- (\gamma_1),\ldots,i_- (\gamma_{2g}), 
z_1 ,\ldots, z_l, 
i_+ (\gamma_1),\ldots,i_+ (\gamma_{2g}) \mid 
r_1, \ldots, r_{2g+l}
\rangle
\end{align}
for some integer $l$. 
That is, it is a finite presentation with deficiency $2g$ 
whose generating set 
contains $i_- (\gamma_1),\ldots,i_- (\gamma_{2g}), 
i_+ (\gamma_1),\ldots,i_+ (\gamma_{2g})$ and is ordered as above. 
Such a presentation always exists (see \cite[Section 4]{gs08}). 
For any admissible presentation, 
define $2g \times (2g+l)$, $l \times (2g+l)$ and 
$2g \times (2g+l)$ matrices $A,B,C$ over $\mathbb{Z} \Gamma$ by 
\[A=\overline{
\sideset{^{\rho}\!}{}
{\mathop{\left({\displaystyle
\frac{\partial r_j}{\partial i_-(\gamma_i)}
} \right)}\nolimits}
}_{\begin{subarray}{c}
{}1 \le i \le 2g\\
1 \le j \le 2g+l
\end{subarray}}, \ \ 
%
B=\overline{
\sideset{^{\rho}\!}{}
{\mathop{\left({\displaystyle
\frac{\partial r_j}{\partial z_i}
} \right)}\nolimits}
}_{\begin{subarray}{c}
{}1 \le i \le l\\
1 \le j \le 2g+l
\end{subarray}}, \ \ 
C=\overline{
\sideset{^{\rho}\!}{}
{\mathop{\left({\displaystyle
\frac{\partial r_j}{\partial i_+(\gamma_i)}
} \right)}\nolimits}
}_{\begin{subarray}{c}
{}1 \le i \le 2g\\
1 \le j \le 2g+l
\end{subarray}}.\]

\begin{proposition}[{\cite[Proposition 4.1]{gs08}}]
\label{prop:MagnusFormula}
As matrices with entries in $\mathcal{K}_\Gamma$, we have:
\begin{itemize}
\item[$(1)$] 
The square matrix $\begin{pmatrix} A \\ B \end{pmatrix}$ 
is invertible and 
$\tau_{\rho}^+ (M_R)=\begin{pmatrix} A \\ B \end{pmatrix}$; and 
\item[$(2)$] $r_\rho(M_R) = 
-C \begin{pmatrix} A \\ B \end{pmatrix}^{-1} \!
\begin{pmatrix} I_{2g} \\ 0_{(l,2g)}\end{pmatrix}$. 
\end{itemize}
\end{proposition}
\begin{remark}\label{rem:Strebel}
We see from Strebel \cite{strebel} 
that for a PTFA group $\Gamma$, 
every matrix with
entries in $\mathbb{Z} \Gamma$ sent to an 
invertible matrix over $\mathbb{Q}$ 
by the augmentation map $\mathbb{Z} \Gamma \to \mathbb{Z}$ 
is invertible over $\mathcal{K}_\Gamma$. 
The first assertion of (1) follows from this fact. 
Indeed, $\begin{pmatrix} A \\ B \end{pmatrix}$ is 
sent to a representation 
matrix of $H_1 (M_R;\mathbb{Q})/
i_+ (H_1 (\Sigma_{g,1};\mathbb{Q}))=0$ 
by the augmentation map, which is invertible over $\mathbb{Q}$. 
\end{remark}
\begin{remark}\label{rem:obstruction}
If $K$ is fibered, 
the complementary sutured manifold for the unique minimal genus Seifert 
surface is a product sutured manifold, 
so that $\Gamma$-torsion is trivial for any $\Gamma$. 
Therefore we can use $\Gamma$-torsion as fibering obstructions of 
homologically fibered knots. Note that 
we can also use the Magnus matrix 
(see \cite[Theorem 4.1]{gs08} and Section \ref{sec:magnus}). 
\end{remark}

By using the above invariants, the factorization formula for 
$\tau_\rho (E(K))$ is given as follows, where 
the statement is simpler than that in \cite{gs08} because we are 
now considering the knot cases only.
\begin{theorem}\label{thm:factorization}
Let $K$ be a rationally homologically fibered knot of genus $g$ and 
let $R$ be a minimal genus Seifert surface of $K$. 
For any non-trivial homomorphism $\rho:G(K) \to \Gamma$ to 
a PTFA group $\Gamma$, a loop $\mu$ representing the meridian of 
$K$ satisfies $\rho(\mu) \neq 1 \in \Gamma 
\subset \mathcal{K}_\Gamma$ and we have a factorization 
\begin{equation}\label{eq2}
\tau_\rho (E(K)) = \frac{\tau_\rho^+ (M_R) \cdot 
(I_{2g} -\rho (\mu) r_\rho (M_R))}{1-\rho(\mu)} \quad 
\in K_1(\mathcal{K}_\Gamma)/\pm\rho(G(K))
\end{equation}
of the torsion $\tau_\rho (E(K))$. 
\end{theorem}

\begin{proof} \ 
First, by passing to the image if necessary, 
we may suppose that $\rho$ is onto. 
This is justified by the facts that any subgroup $\Gamma'$ 
of a PTFA group $\Gamma$ is again PTFA and that the torsion 
is invariant under the field extension $\mathcal{K}_{\Gamma'} 
\hookrightarrow \mathcal{K}_\Gamma$. 

By the definition of PTFA groups, we see that there exists a surjective 
homomorphism $\Gamma \to \mathbb{Z}$. Then the composite 
$G(K) \xrightarrow{\rho} \Gamma \to \mathbb{Z}$ is also 
surjective and it coincides with the abelianization map 
of $G(K)$ up to sign. Hence $\rho (\mu) \neq 1 \in \Gamma$. 

The rest of the proof is almost identical to the argument in 
\cite[Section 5]{gs08} (see also the argument of 
Friedl \cite[Section 6]{fri}). For convenience, 
we repeat it here in a simplified form. 

Given an admissible presentation of $\pi_1 (M_R)$ 
as in (\ref{admissible}), we denote it briefly by 
\[\pi_1 (M_R) \cong \langle i_- (\overrightarrow{\gamma}), 
\overrightarrow{z}, i_+ (\overrightarrow{\gamma}) \mid 
\overrightarrow{r} \rangle.\]
A usual computation gives 
\[G(K) \cong \langle i_- (\overrightarrow{\gamma}),
\overrightarrow{z},i_+ (\overrightarrow{\gamma}),\mu \mid 
\overrightarrow{r}, i_- (\overrightarrow{\gamma}) \, \mu \, 
i_+(\overrightarrow{\gamma})^{-1} \mu^{-1} \rangle.\]
From this presentation, we construct a 2-complex $X(K)$ 
consisting of one 0-cell, 
one 1-cell for each generator and one 2-cell 
for each relation with an attaching map according to the word. 
We can check that $E(K)$ and $X(K)$ 
are simple homotopy equivalent 
(see \cite[Lemma 5.1]{gs08}). 

The $\mathcal{K}_\Gamma$-rank 
of $C_i (X(K);\mathcal{K}_\Gamma)$ 
and the $\mathbb{Z}$-rank $C_i (X(K))$ are the same and 
their degree $0, 1, 2$ parts are 
given by $1, 4g+l+1, 4g+l$. 
The map $\partial_2: C_2 (X(K);\mathcal{K}_\Gamma) 
\cong \mathcal{K}_\Gamma^{4g+l} 
\to C_1 (X(K);\mathcal{K}_\Gamma) \cong \mathcal{K}_\Gamma^{4g+l+1}$ 
is represented by the matrix 
\[D_2:=
\begin{pmatrix}
A & I_{2g} \\
B & 0_{(l,2g)} \\
C & -\rho (\mu)^{-1} I_{2g} \\
0_{(1,2g+l)} & \ast \ \ast \ \cdots \ \ast
\end{pmatrix}. \]
Consider 
the matrix $D_2^\mu$ obtained from 
$D_2$ by deleting the last row. 
By fundamental transformations of matrices, we have 
\begin{align*}
D_2^\mu &=
\begin{pmatrix}
A \ & I_{2g} \\
B \ & 0_{(l,2g)}\\
C \ & -\rho (\mu)^{-1} I_{2g}
\end{pmatrix}
\to 
\begin{pmatrix}
A + \rho (\mu) C & 0_{2g} \\
B & 0_{(l,2g)}\\
C & -\rho (\mu)^{-1} I_{2g}
\end{pmatrix}\\
&\to 
\begin{pmatrix}
A + \rho (\mu) C & 0_{2g} \\ 
B & 0_{(l,2g)} \\
0_{(2g,2g+l)} & -\rho (\mu)^{-1} I_{2g}
\end{pmatrix} =:D
\end{align*}
Here the above matrices are with entries in $\mathbb{Z} \Gamma$ and 
we apply the augmentation map $\mathbb{Z} \Gamma \to \mathbb{Z}$ to 
$D$. Then we have a matrix representing 
$\partial_2: C_2 (X(K)) \to C_1 (X(K))/ \langle \mu \rangle$, 
which can be easily seen to be 
invertible over $\mathbb{Q}$. 
Hence $D$ (and also $D_2^\mu$) is invertible over $\mathcal{K}_\Gamma$ 
as mentioned in Remark \ref{rem:Strebel}. 


Now we compute $\tau_\rho (E(K))=\tau (C_\ast (X(K);\mathcal{K}_\Gamma))$. 
By the cell structure of $X(K)$, 
\[\tau_\rho (E(K)) = D_2^\mu \cdot (1-\rho (\mu)^{-1})^{-1}\] 
holds. Then as elements 
in $K_1(\mathcal{K}_\Gamma)/\pm\rho(G(K))$, we have 
\begin{align*}
D_2^\mu &= D =\begin{pmatrix}
A + \rho (\mu) C \\ B 
\end{pmatrix}=\begin{pmatrix}
I_{2g} - \rho (\mu) r_{\rho} (M_R) & 
\ -\rho (\mu) Z \\
0_{(l,2g)} & \ I_l
\end{pmatrix}
\begin{pmatrix}
A \\ B
\end{pmatrix}\\
&=(I_{2g} - \rho (\mu) r_{\rho} (M_R)) 
\begin{pmatrix}
A \\ B 
\end{pmatrix},
\end{align*}
\noindent
where we used 
\[\begin{pmatrix}
A + \rho (\mu) C \\ B 
\end{pmatrix}
=\begin{pmatrix}
A \\ B
\end{pmatrix}
-\rho (\mu)
\begin{pmatrix}
r_{\rho} (M_R) \quad Z \\
0_{(l,2g+l)}
\end{pmatrix}
\begin{pmatrix}
A \\ B
\end{pmatrix}\]
at the third equality and $Z$ is defined by the formula 
$(r_\rho (M_R) \quad Z)= 
-C \begin{pmatrix} A \\ B \end{pmatrix}^{-1}$ 
(see Proposition \ref{prop:MagnusFormula} (2)). 
This completes the proof. 
\end{proof}

\noindent
When we take the abelianization map 
$\rho_1: G(K) 
\to \langle t \rangle \subset \mathbb{Q} (t)$ 
as $\rho$, the formula (\ref{eq1}) is recovered. 

\begin{remark}
Factorizations of (higher-order) Alexander invariants 
into some torsions and ``monodromy'' information 
appear in various contexts such as Morse-Novikov theory and 
the theory of string links. 
For example, see Hutchings-Lee \cite{hl,hl2}, 
Goda-Matsuda-Pajitnov \cite{GMP}, 
Kitayama \cite{kitayama} and Kirk-Livingston-Wang \cite{klw}. 
It would be interesting to compare these factorization formulas 
in an appropriate situation. 
\end{remark}


\section{A sample calculation }\label{sec:sample}

Although all the ingredients in the formula $(\ref{eq2})$ 
are determined by information on fundamental groups, 
it is difficult to compute them explicitly 
because of the non-commutativity of $\mathcal{K}_\Gamma$ except 
in some special cases including the following. 

Let $K$ be a homologically fibered knot with 
a minimal genus Seifert surface $R$ and 
let $M_{R}$ be the sutured manifold for $R$. 
Consider the group extension 
\begin{equation}\label{eq:seq}
1\longrightarrow G(K)'/G(K)''
 \longrightarrow D_2 (K) 
 \longrightarrow G(K)/G(K)'=H_1(E(K))\cong \mathbb{Z}
 \longrightarrow 1
\end{equation}
relating to the metabelian quotient $D_2 (K) :=G(K)/G(K)''$ 
of $G(K)$. 
We have 
\[G(K)'/G(K)'' \cong H_{1}(R) \cong H_{1}(M_{R})\]
since it coincides with 
the first homology of the infinite cyclic 
covering of $E(K)$, which can be seen as 
the product (as homology cylinders) 
of infinitely many copies of $M_R$. 
Let $\rho_2$ be the natural projection 
\[\rho_2: G(K) \longrightarrow D_2 (K). \]
It is known that $D_2 (K)$ is PTFA (see Strebel \cite{strebel}), 
so that $\mathcal{K}_{D_2 (K)}$ is defined. 
Then, it follows from the Proposition \ref{prop:MagnusFormula} that 
$\tau_{\rho_2}^+ (M_R)$ and $r_{\rho_2} (M_{R})$
can be computed by calculations on the {\it commutative} subfield 
$\mathcal{K}_{H_1 (M_R)}$ of $\mathcal{K}_{D_2 (K)}$, 
and therefore we can carry it out. 

\begin{figure}[htbp]
\begin{center}
\includegraphics[width=0.8\textwidth]{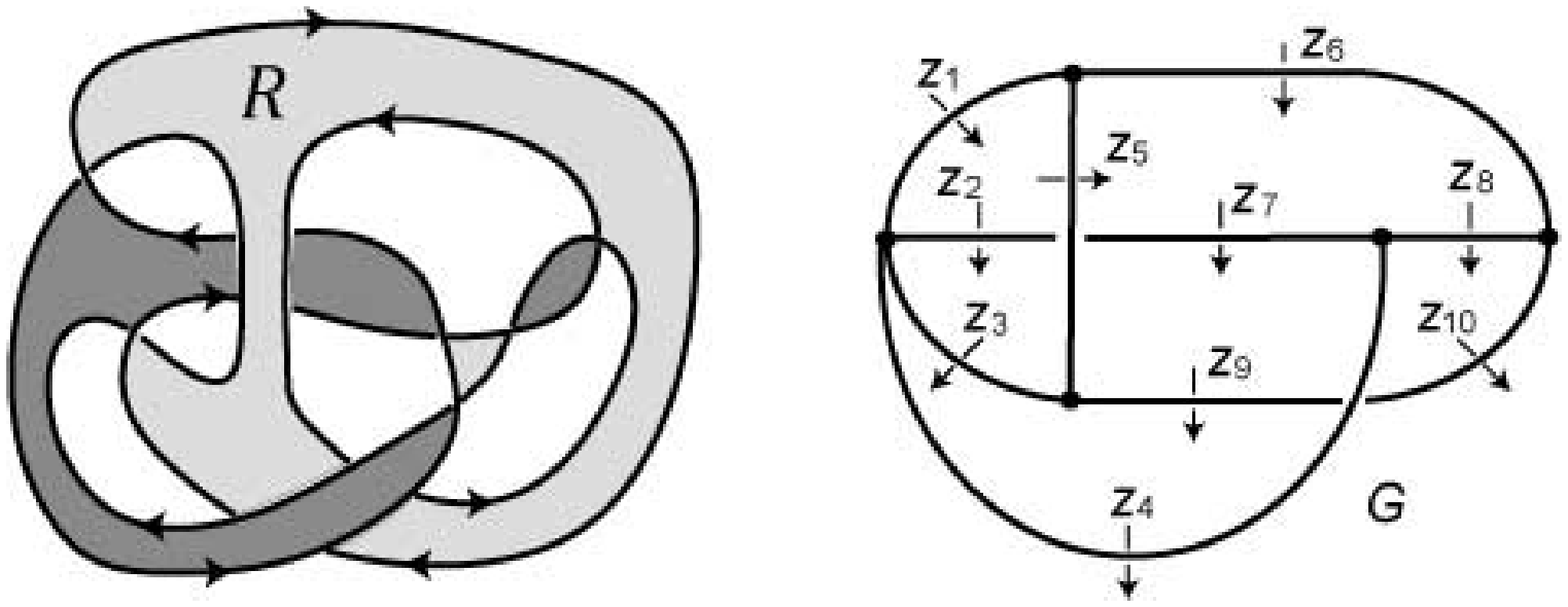}
\end{center}
\caption{}
\label{fig:0057knot}
\end{figure}

Let us see an example of calculations of our invariants. 
Let $K$ be the knot obtained as the boundary of the Seifert surface $R$ 
illustrated in Figure \ref{fig:0057knot}. 
We can easily compute that $\Delta_K (t)=1-2t+3t^2-2t^3+t^4$ 
and the genus of $R$ is $2$. 
Hence $K$ is a homologically fibered knot and 
$R$ is of minimal genus. 
The graph $G$ in right hand side of Figure \ref{fig:0057knot} 
is obtained from $R$ by a deformation retract. 
Thus $\pi_{1}(M_R)\cong \pi_{1}(S^3-\overset{\circ}{N}(G))$. 
Then $\pi_1(M_{R})$ has a presentation:
\[{\small \left\langle 
\begin{array}{c|l}
z_{1},z_{2}, \ldots, z_{10} \ & \ 
\begin{array}{ll}
z_{1}z_{5}z_{6}^{-1},\,z_{2}z_{3}z_{4}z_{1},\,z_{3}z_{9}^{-1}z_{5}^{-1},\,
z_{7}z_{4}z_{8}^{-1},\\
z_{8}z_{10}z_{6},\,
z_{2}z_{5}z_{7}^{-1}z_{5}^{-1},\,z_{9}z_{4}z_{10}^{-1}z_{4}^{-1}
\end{array}
\end{array}
\right\rangle.
}\]
We can drop the last relation $z_{9}z_{4}z_{10}^{-1}z_{4}^{-1}$ 
because it is derived from the others.

\begin{figure}[h]
\begin{center}
\includegraphics[width=0.4\textwidth]{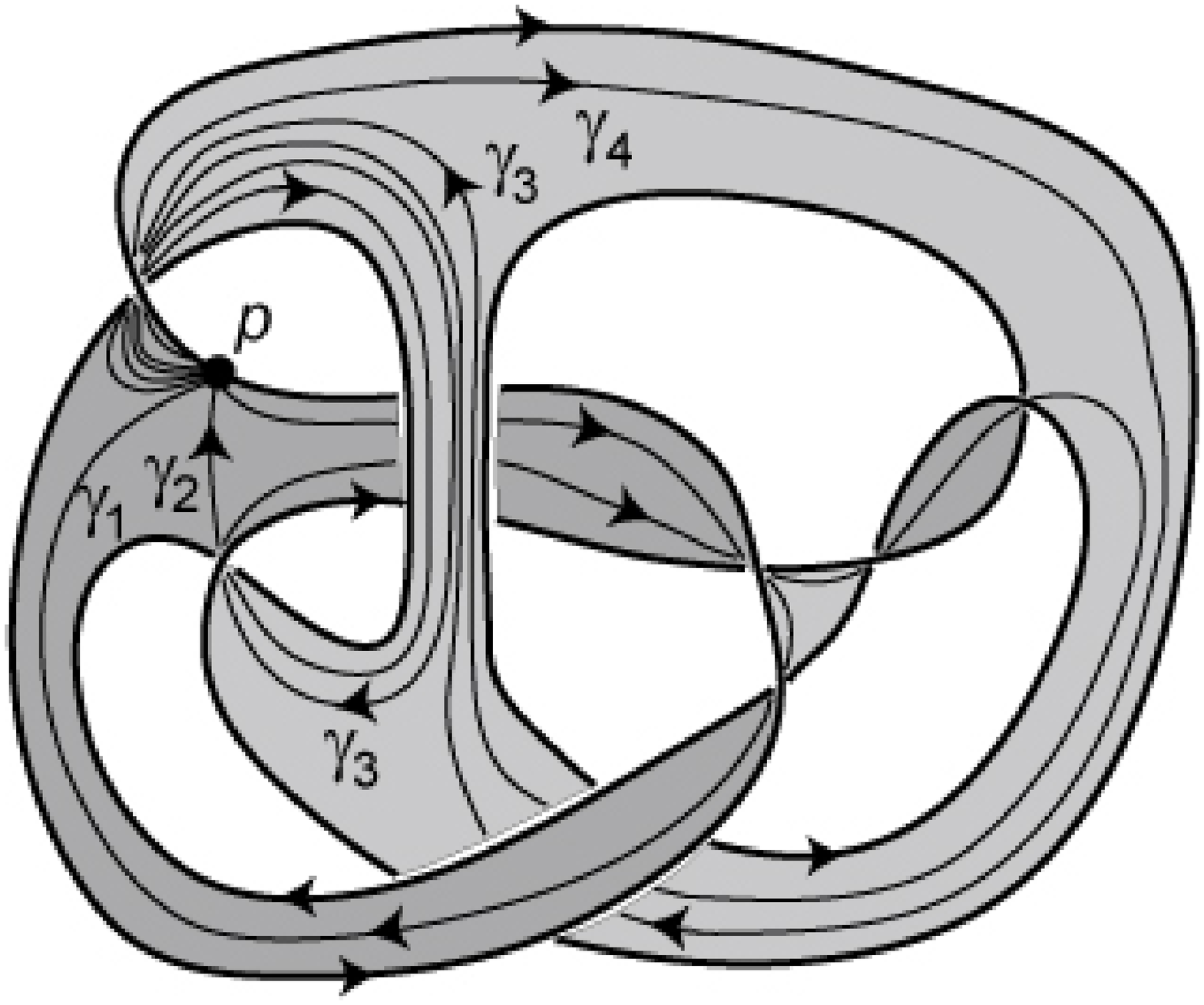}
\end{center}
\caption{}
\label{fig:basis}
\end{figure}

We take a spine of $R$ as in Figure \ref{fig:basis}, 
by which we can fix an identification of $\Sigma_{g,1}$ and $R$. 
A direct computation shows that 
\begin{align*}
i_-(\gamma_1)&=z_{5}z_{1}, & 
i_-(\gamma_2)&=z_{2}^{-1}, &
i_-(\gamma_3)&=z_{5}z_{7}^{-1}z_{8}^{-1}z_{4}^{-1}, & 
i_-(\gamma_4)&=z_{4}^{-1}, \\
i_+(\gamma_1)&=z_{5}, & 
i_+(\gamma_2)&=z_{6}z_{9}, &
i_+(\gamma_3)&=z_{6}z_{5}^{-1}z_{3}z_{5}z_{7}^{-1}z_{4}^{-1}z_{6}^{-1}, & 
i_+(\gamma_4)&=z_{6}z_{7}z_{6}^{-1}. 
\end{align*}
\noindent
Here the darker color in $R$ is the $+$-side.
Then, we obtain an admissible presentation of $\pi_1 (M_{R})$:
{\small \begin{center}
\begin{tabular}{lcr}
\multicolumn{1}{c}{Generators} & 
$i_{-}(\gamma_1),\ldots, i_{-}(\gamma_4),\, 
z_{1},\ldots ,z_{10},\,
i_{+}(\gamma_{1}),\ldots,i_{+}(\gamma_{4})$ \\ 
\multicolumn{1}{c}{Relations} & 
$z_{1}z_{5}z_{6}^{-1},\,z_{2}z_{3}z_{4}z_{1},\,z_{3}z_{9}^{-1}z_{5}^{-1},\,
z_{7}z_{4}z_{8}^{-1},\,z_{8}z_{10}z_{6},\,
z_{2}z_{5}z_{7}^{-1}z_{5}^{-1},
$ \\
	&
$i_{-}(\gamma_{1})z_{1}^{-1}z_{5}^{-1},\, i_{-}(\gamma_{2})z_{2},\,
i_{-}(\gamma_{3})z_{4}z_{8}z_{7}z_{5}^{-1},\, i_{-}(\gamma_{4})z_{4},$ \\
	&
$i_{+}(\gamma_{1})z_{5}^{-1}, \,i_{+}(\gamma_{2})z_{9}^{-1}z_{6}^{-1},\,
 i_{+}(\gamma_{3})z_{6}z_{4}z_{7}z_{5}^{-1}z_{3}^{-1}z_{5}z_{6}^{-1},\,
 i_{+}(\gamma_{4})z_{6}z_{7}^{-1}z_{6}^{-1}$
\end{tabular}
\end{center}
}

\begin{figure}[htbp]
\begin{center}
\includegraphics[width=0.65\textwidth]{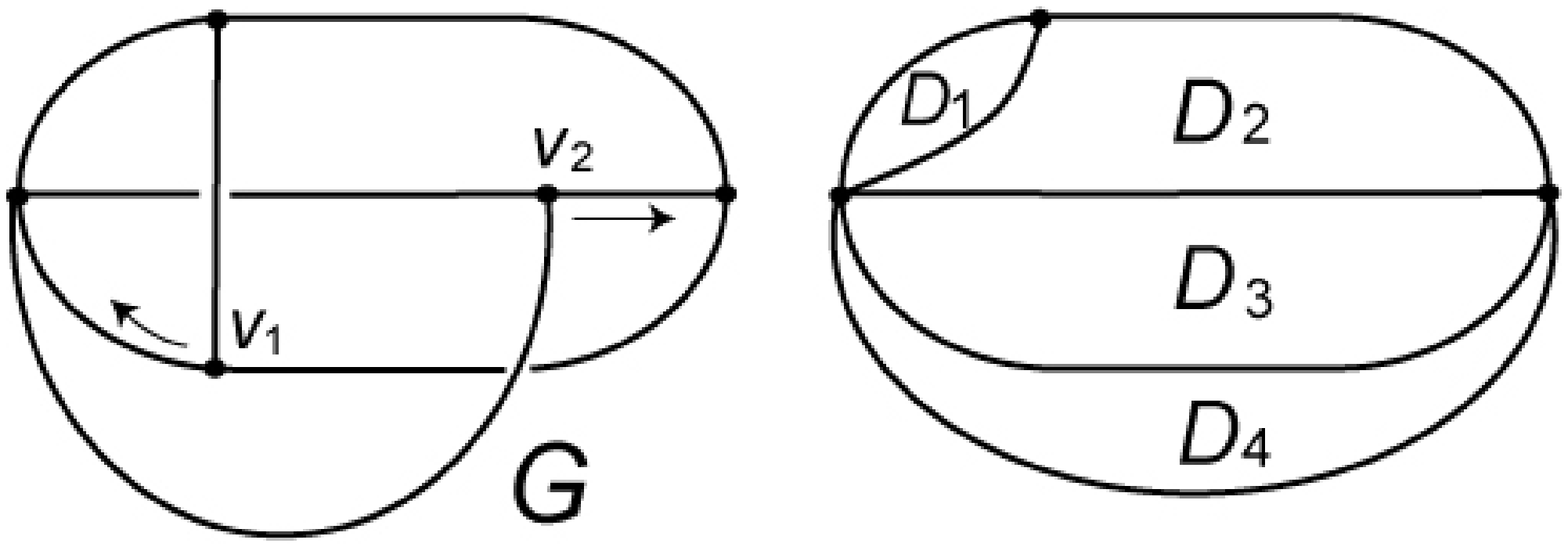}
\end{center}
\caption{}
\label{fig:homology}
\end{figure}

By sliding the edges $v_{1}$ and $v_{2}$ of $G$ as in 
Figure \ref{fig:homology}, we obtain a graph 
whose complement is a genus 4 handlebody. 
This means that the complement of $G$ (and hence $M_R$) 
is homeomorphic to a genus 4 handlebody. 

%

Let $D_1,\ldots , D_4$ be the meridian disks of the handlebody 
as illustrated in the figure. Put 
$x_1:=z_{1}^{-1}$, $x_{2}=z_{6}^{-1}$, $x_{3}:=(z_{6}z_{7})^{-1}$ 
and $x_{4}:=z_{4}$, 
where $x_i$ is a loop intersecting $D_i$ transversely 
in one point from the above to the down side 
in Figure \ref{fig:homology} and
is disjoint from $D_{j}\,(i\neq j)$. 
By using them, 
we have the following simplified 
admissible presentation of $\pi_1(M_{R})$: 
{\small \begin{center}
\begin{tabular}{lcr}
\multicolumn{1}{c}{Generators} & 
$i_{-}(\gamma_1),\ldots, i_{-}(\gamma_4),\, 
x_{1},x_{2},x_{3} ,x_{4},\,
i_{+}(\gamma_{1}),\ldots,i_{+}(\gamma_{4})$ \\ 
\multicolumn{1}{c}{Relations} & 
$i_{-}(\gamma_{1})x_{1}x_{2}x_{1}^{-1}$, \,
$i_{-}(\gamma_{2})x_{1}x_{3}^{-1}x_{2}x_{1}^{-1}$, \,
$i_{-}(\gamma_{3})x_{4}x_{2}x_{3}^{-1}x_{4}x_{2}x_{3}^{-1}x_{2}x_{1}^{-1}$,\\ 
 & $i_{-}(\gamma_{4})x_{4}$, \,$i_{+}(\gamma_{1})x_{2}x_{1}^{-1},$ \,
$i_{+}(\gamma_{2})x_{4}x_{3}^{-1}x_{2}$,\\
 & $i_{+}(\gamma_{3})x_{2}^{-1}x_{4}x_{2}x_{3}^{-1}
 x_{2}x_{1}^{-1}x_{4}x_{3}^{-1}x_{2}$, \,$i_{+}(\gamma_{4})x_{2}^{-1}x_{3}$
\end{tabular}
\end{center}
}
We write $r_1,\ldots ,r_8$ for these relations in order. 
Recall that 
$\mathcal K_{H_1 (M_R)}$ is isomorphic to the field of 
rational functions with variables $x_1,\ldots , x_4$, where 
we use the same notation for the image of $x_i$ by the 
abelianization map $\pi_1 (M_R) \to H_1 (M_R)$. 
Then we have
\noindent
\[\begin{pmatrix}A \\ B\\ C\end{pmatrix}=
%
\begin{pmatrix}
I_4 & 0_4 \\
G_1 & G_2 \\ 
0_4 & I_4
\end{pmatrix},\] 
where $G_1 = \begin{pmatrix}
g_{11}&g_{12}&g_{13}&g_{14}\\
g_{21}&g_{22}&g_{23}&g_{24}\\
g_{31}&g_{32}&g_{33}&g_{34}\\
g_{41}&g_{42}&g_{43}&g_{44}
\end{pmatrix}$ and 
$G_2 = \begin{pmatrix}
g_{15}&g_{16}&g_{17}&g_{18}\\
g_{25}&g_{26}&g_{27}&g_{28}\\
g_{35}&g_{36}&g_{37}&g_{38}\\
g_{45}&g_{46}&g_{47}&g_{48}
\end{pmatrix}$ with 
$g_{ij}=
\rho\left(
\overline{\displaystyle\frac{\partial r_{j}}{\partial x_{i}}}
\right)$. 
Thus 
\[\tau_{\rho_2}^+ (M_{R}) = \begin{pmatrix}A \\ B\end{pmatrix}
=
\begin{pmatrix}
I_4 & 0_4 \\
G_1 & G_2
\end{pmatrix}.\]
As a torsion, 
it is equivalent to $G_2$, 
where 
\begin{alignat*}{3}
   g_{15} &= -1,
& \quad g_{16} &= 0,
& \quad g_{18} &= 0,\\ 
   g_{25} &= x_{1}^{-1}x_{2},
& \quad g_{26} &= x_{2},
& \quad g_{28} &= -x_{3},\\
   g_{35} &= 0, 
& \quad g_{36} &= -x_{2},
& \quad g_{38} &= x_{3},\\ 
   g_{45} &= 0,
& \quad g_{46} &= x_{2}x_{3}^{-1}x_{4}, 
& \quad g_{48} &= 0,
\end{alignat*}
\begin{align*}
g_{17} &= -x_{2}x_{3}^{-1}x_{4}, \qquad 
g_{27} = x_{2}+x_{1}^{-1}x_{2}^{2}x_{3}^{-1}x_{4}+x_{1}^{-1}x_{2}^{3}x_{3}^{-2}x_{4}-x_{1}^{-1}x_{2}^{3}x_{3}^{-2}x_{4}^{2},\\
g_{37} &= -x_{2}-x_{1}^{-1}x_{2}^{2}x_{3}^{-1}x_{4}, \qquad 
g_{47} = x_{2}x_{3}^{-1}x_{4}+x_{1}^{-1}x_{2}^{3}x_{3}^{-2}x_{4}^{2}.
\end{align*}
Then we have
\[\det (\tau_{\rho_2}^+ (M_{R}))=
\det (G_2)
=
-\frac{x_{2}^{3}x_{4}^{2}}{x_{1}x_{3}^{2}}(x_{2}-x_{3}-x_{2}x_{4}).
\]
The Magnus matrix $r_{\rho_2} (M_R)$ can be computed by the formula in 
Proposition \ref{prop:MagnusFormula} (2). 
However we omit it here. 

\begin{remark}
From an admissible presentation, 
we can use the Mathematica program given 
in Section \ref{sec:program} for
calculations of $\tau_{\rho_2}^+ (M_{R})$ and $r_{\rho_2} (M_{R})$. 
Note that the program uses 
$\{i_+(\gamma_1), i_+(\gamma_2),\ldots,i_+(\gamma_{2g})\}$ as 
a basis of $H_1 (M_R)$. In the above example, 
\[x_{1}=\gamma_{2}^{-2}\gamma_{3}, \quad 
x_{2}=\gamma_{1}^{-1}\gamma_{2}^{-2}\gamma_{3}, \quad 
x_{3}=\gamma_{1}^{-1}\gamma_{2}^{-2}\gamma_{3}\gamma_{4}^{-1}, \quad 
x_{4}=\gamma_{2}^{-1}\gamma_{4}^{-1},\]
where $\gamma_j$ denotes $i_+ (\gamma_j)$, and we have 
$\det (\tau_{\rho_2}^+ (M_{R}))=
\displaystyle\frac{\gamma_{3}}{\gamma_{1}^2 \gamma_{2}^5 \gamma_{4}}
(1+\gamma_{2}-\gamma_{2}\gamma_{4})$. 
\end{remark}


\section{Homologically fibered knots with 12-crossings}\label{sec:HFK12}

\begin{figure}[htbp]
\begin{center}
\includegraphics[width=0.99\textwidth]{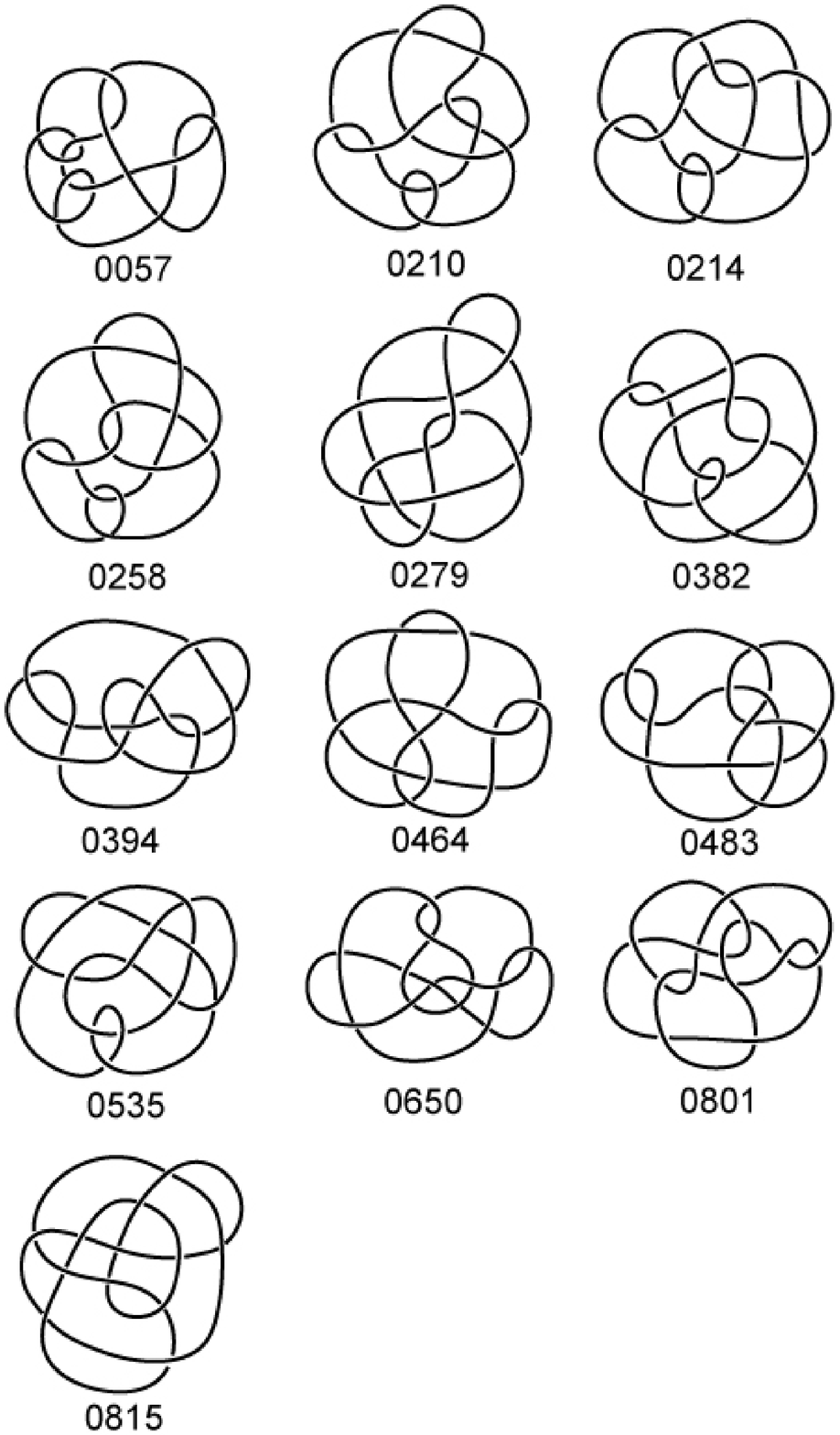}
\end{center}
\caption{12-crossings non-fibered homologically fibered knots}
\label{fig:all}
\end{figure}

It is known that all homologically fibered knots are fibered
among prime knots with at most 11-crossings. 
On the other hand, Friedl-Kim \cite{fk} showed that 
there are 13 non-fibered homologically 
fibered knots with 12-crossings 
by using the twisted Alexander invariant.
See Figure \ref{fig:all} and Table \ref{table:nonfibered}.  
In this section, we list admissible presentations and 
the torsion $\tau^{+}_{\rho_2}$ for sutured manifolds 
associated with minimal genus Seifert surfaces 
illustrated in Figure $\ref{fig:0210SG}, \ldots, \ref{fig:0815SG}$.  
As a by-product, we observe that $\tau^{+}_{\rho_2}$ 
can also detect the non-fiberedness of all these knots. 
In the forthcoming paper \cite{gs09}, 
we will obtain the same result by using 
Johnson homomorphisms as a fibering obstruction.

It is easy to see that the complements of the Seifert 
surfaces for knots 0210, 0214, 0382 and 0394 are handlebodies. 
(In \cite{cl}, a non-alternating prime knot with 12-crossings 
is denoted by $12n_{-}P$.
We refer only the number $P$ in this section.)
Hence, we take free generators corresponding to disks $z_{i}$ 
in each figure,  
which run from the upside to the downside of the diagrams. 
As for the other knots, we have the admissible presentations 
by the same method as in Section \ref{sec:sample}. 

\begin{table}[htbp]
\begin{tabular}{|c|c|c|}
\hline
Knot & Genus & Alexander polynomial \\
\hline\hline
0057 & 2 & $1-2t+3t^2-2t^3+t^4$\\
\hline
0210, 0214 & 3 & $1-t-t^2+3t^3-t^4-t^5+t^6$\\
\hline
0258, 0464, 0483 & 2 & $1-4t+5t^2-4t^3+t^4$ \\
\hline 
0279, 0394 & 2 & $1-6t+11t^2-6t^3+t^4$\\
\hline 
0382, 0801 & 2 & $1-5t+7t^2-5t^3+t^4$\\
\hline 
0535 & 2 & $1-7t+11t^2-7t^3+t^4$\\
\hline 
0650 & 2 & $1-4t+7t^2-4t^4+t^4$\\
\hline 
0815 & 2 & $1-2t+t^2-2t^3+t^4$\\
\hline
\end{tabular}
\caption{Non-fibered homologically fibered knots with 12-crossings}
\label{table:nonfibered} 
\end{table}

The following are admissible presentations and the determinant of 
the torsion $\tau^{+}_{\rho_2} (M_R)$, where we use 
$\{i_+(\gamma_1), i_+(\gamma_2),\ldots, i_+(\gamma_{2g})\}$ 
as a basis of 
$H_1 (M_R)$ and 
denote $i_+(\gamma_j)$ by $\gamma_j$ for simplicity. 
Note that the example in Section \ref{sec:sample} 
is about the knot 0057, and we omit it here. 


\begin{center}
\begin{tabular}{lcr}
\multicolumn{2}{c}{0210}\\
\multicolumn{1}{c}{Generators} & 
$i_{-}(\gamma_{1}),\ldots, i_{-}(\gamma_{6}),\, 
 z_{1},\ldots z_{6},\,
 i_{+}(\gamma_{1}),\ldots,i_{+}(\gamma_{6})$ \\ 
\multicolumn{1}{c}{Relations} & 
$
i_{-}(\gamma_{1})z_{3}^{-1}z_{4},\, 
i_{-}(\gamma_{2})z_{3}^{-2}z_{2},\,
i_{-}(\gamma_{3})z_{5}^{-1}z_{3}^{-1}z_{2},\,
i_{-}(\gamma_{4})z_{2}^{-1}z_{1}z_{6}^{-1}z_{5}z_{6}^{-1}z_{5},\,
$ \\
	& 
$
i_{-}(\gamma_{5})z_{5}^{-1}z_{6}z_{5}^{-1}z_{1}z_{6}^{-1}z_{5}z_{6}^{-1}z_{5},\,
i_{-}(\gamma_{6})z_{5}^{-1}z_{6}z_{5}^{-1}z_{1}z_{3}^{-1}z_{5}z_{6}^{-1}z_{5},\,
$ \\
	&
$
i_{+}(\gamma_{1})z_{4}, \,
i_{+}(\gamma_{2})z_{4}z_{3}^{-1}z_{2}z_{3}^{-1},\,
i_{+}(\gamma_{3})z_{6}^{-1}z_{2}z_{3}^{-1},\,
i_{+}(\gamma_{4})z_{5}z_{2}^{-1}z_{1}z_{6}^{-1}z_{5},\,
$ \\
	& 
$
i_{+}(\gamma_{5})z_{5}^{-1}z_{6}z_{2}^{-1}z_{1}z_{6}^{-1}z_{5},\,
i_{+}(\gamma_{6})z_{5}^{-1}z_{6}z_{3}^{-1}z_{5}z_{6}^{-1}z_{5}\,
$
\\
\multicolumn{1}{c}{Torsion $\tau^{+}_{\rho_2}$} & 
$\displaystyle{-\frac{\gamma_1^5\gamma_3^3\gamma_5^4\gamma_6^7}{\gamma_2^6\gamma_4^6}
+\frac{\gamma_1^6\gamma_3^4\gamma_5^4\gamma_6^7}{\gamma_2^7\gamma_4^6}
-\frac{\gamma_1^6\gamma_3^4\gamma_5^4\gamma_6^8}{\gamma_2^7\gamma_4^6}}$\\
\end{tabular}
\begin{figure*}[h]
\includegraphics[width=0.99\textwidth]{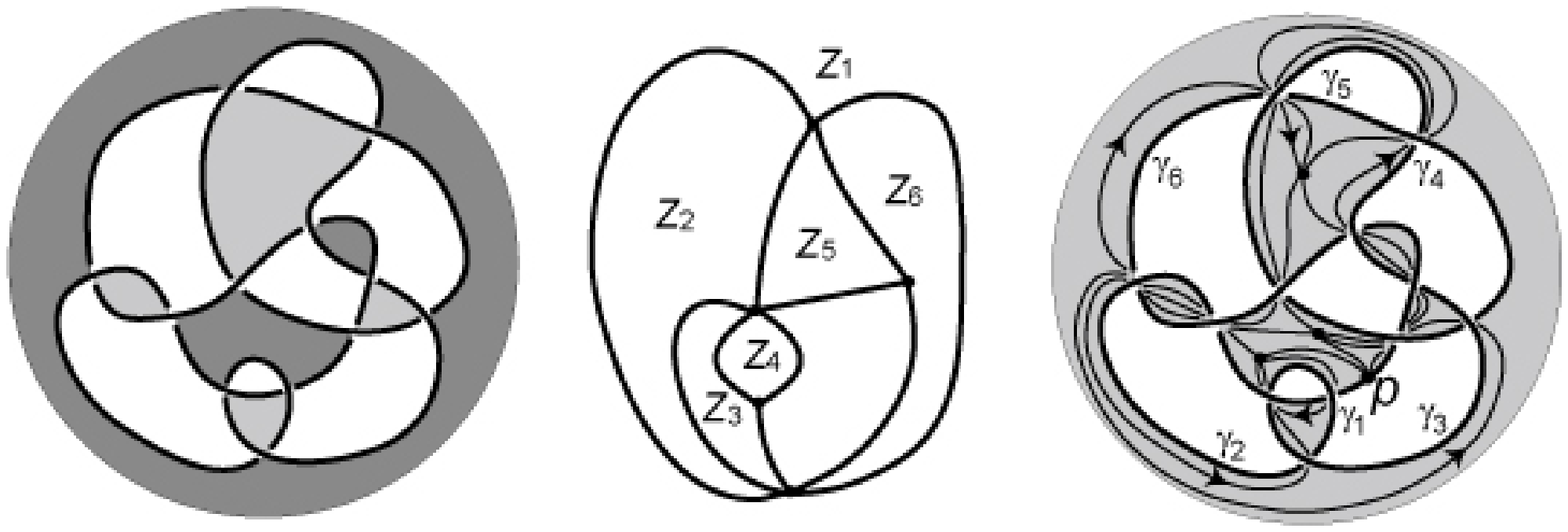}
\caption{0210}
\label{fig:0210SG}
\end{figure*}
\end{center}

\begin{center}
\begin{tabular}{lcr}
\multicolumn{2}{c}{0214}\\
\multicolumn{1}{c}{Generators} & 
$i_{-}(\gamma_{1}),\ldots, i_{-}(\gamma_{6}),\, 
 z_{1},\ldots z_{6},\,
 i_{+}(\gamma_{1}),\ldots,i_{+}(\gamma_{6})$ \\ 
\multicolumn{1}{c}{Relations} & 
$
i_{-}(\gamma_{1})z_{2}z_{3}^{-1}z_{2}^{-1},\, 
i_{-}(\gamma_{2})z_{2}z_{1}^{-1}z_{2},\,
i_{-}(\gamma_{3})z_{5}^{-1}z_{1}^{-1}z_{2},\,
$ \\
	& 
$
i_{-}(\gamma_{4})z_{6}^{-1}z_{1}z_{3}^{-1}z_{5},\,
i_{-}(\gamma_{5})z_{5}^{-1}z_{4}z_{3}^{-1}z_{1}z_{3}^{-1}z_{5},\,
i_{-}(\gamma_{6})z_{5}^{-1}z_{4},\,
$ \\
	&
$
i_{+}(\gamma_{1})z_{2}^{2}z_{3}^{-1}z_{2}^{-1}, \,
i_{+}(\gamma_{2})z_{2}^{2}z_{6}^{-1},\,
i_{+}(\gamma_{3})z_{1}^{-1}z_{2}z_{6}^{-1},\,
$ \\
	& 
$
i_{+}(\gamma_{4})z_{5}z_{3}^{-1}z_{1},\,
i_{+}(\gamma_{5})z_{3}^{-1}z_{5}z_{3}^{-1}z_{1},\,
i_{+}(\gamma_{6})z_{3}^{-1}z_{4}\,
$
\\
\multicolumn{1}{c}{Torsion $\tau^{+}_{\rho_2}$} & 
$\displaystyle{\frac{1}{\gamma_2\gamma_4^2\gamma_6}
-\frac{\gamma_1}{\gamma_2\gamma_4^2\gamma_6}
+\frac{\gamma_1}{\gamma_2\gamma_4\gamma_5\gamma_6}}$ \\
\end{tabular}
\begin{figure*}[h]
\includegraphics[width=0.99\textwidth]{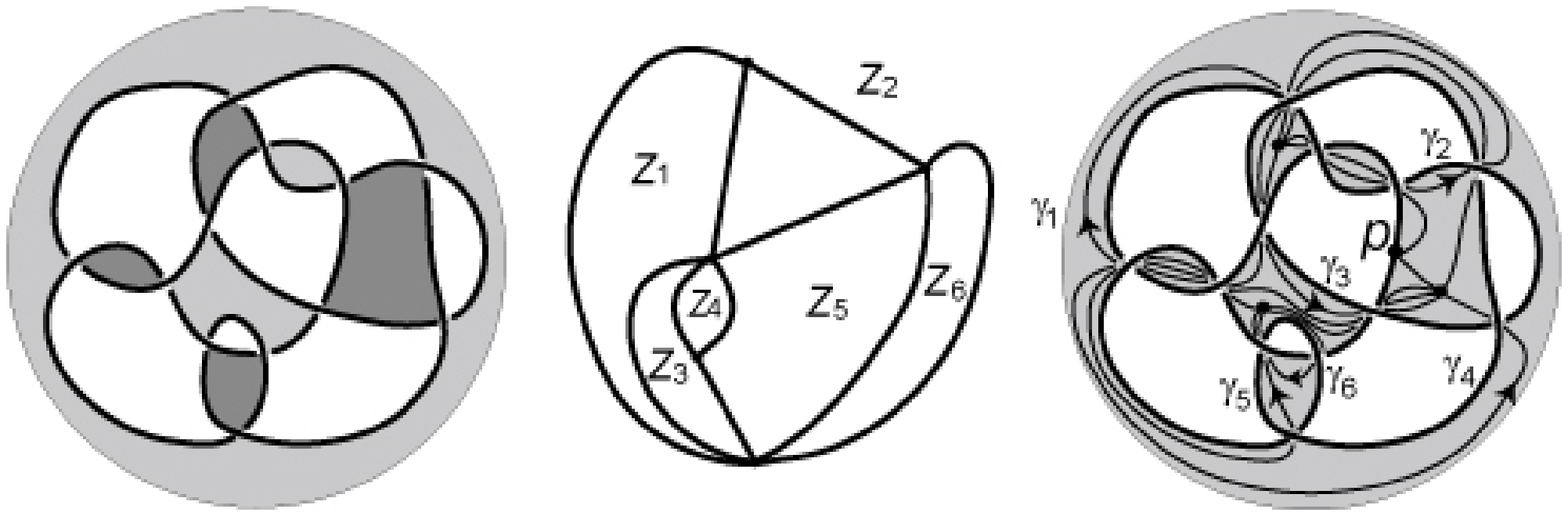}
\caption{0214}
\label{fig:0214SG}
\end{figure*}
\end{center}

\begin{center}
\begin{tabular}{lcr}
\multicolumn{2}{c}{0258}\\
\multicolumn{1}{c}{Generators} & 
$i_{-}(\gamma_1),\ldots, i_{-}(\gamma_4),\, 
z_{1},\ldots ,z_{7},\,
i_{+}(\gamma_{1}),\ldots,i_{+}(\gamma_{4})$ \\ 
\multicolumn{1}{c}{Relations} & 
$z_{1}z_{2}z_{3}z_{4}$, $z_{1}z_{2}z_{4}z_{6}^{-1}z_{7}^{-1}$, 
$z_{7}z_{6}z_{5}$, $i_{-}(\gamma_{1})z_{7}z_{6}z_{7}^{-1}$, \\ 
& $i_{-}(\gamma_{2})z_{7}z_{6}z_{5}^{-1}z_{4}z_{6}^{-1}z_{7}^{-1}$, 
$i_{-}(\gamma_{3})z_{1}z_{2}^{2}z_{4}^{2}z_{6}^{-1}z_{7}^{-1}$, 
$i_{-}(\gamma_{4})z_{1}z_{2}^{2}z_{1}^{-2},$ \\
 &
$i_{+}(\gamma_{1})z_{7}^{-1}, \,i_{+}(\gamma_{2})z_{6}z_{4},\,
 i_{+}(\gamma_{3})z_{2}z_{1}^{-1}z_{4},\,
 i_{+}(\gamma_{4})z_{2}z_{1}^{-2}\,$
\\
\multicolumn{1}{c}{Torsion $\tau^{+}_{\rho_2}$} & 
$\displaystyle{-\frac{\gamma_2^{5}\gamma_4^{7}}{\gamma_1^{6}\gamma_3^{12}}
+\frac{\gamma_2^{6}\gamma_4^{8}}{\gamma_1^{7}\gamma_3^{13}}
-\frac{\gamma_2^{6}\gamma_4^{9}}{\gamma_1^{7}\gamma_3^{14}}}$\\
\end{tabular}
\begin{figure*}[h]
\includegraphics[width=0.99\textwidth]{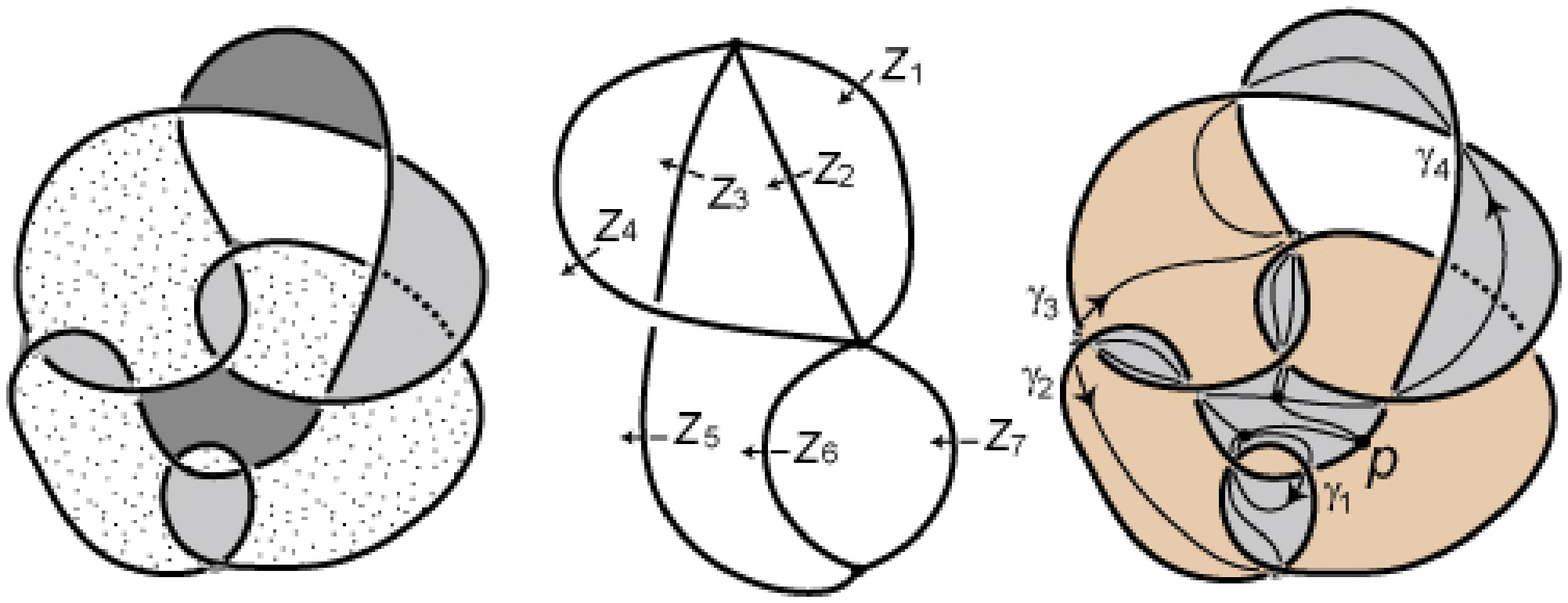}
\caption{0258}
\label{fig:0258SG}
\end{figure*}
\end{center}

\begin{center}
\begin{tabular}{lcr}
\multicolumn{2}{c}{0279}\\
\multicolumn{1}{c}{Generators} & 
$i_{-}(\gamma_1),\ldots, i_{-}(\gamma_4),\, 
z_{1},\ldots ,z_{9},\,
i_{+}(\gamma_{1}),\ldots,i_{+}(\gamma_{4})$ \\ 
\multicolumn{1}{c}{Relations} & 
$
 z_{1}z_{2}z_{4},\,z_{1}z_{3}^{-1}z_{2}z_{9}^{-1},\,z_{5}z_{8}^{-1}z_{6}^{-1},\,
 z_{6}z_{7}z_{8}z_{9},\,z_{2}^{-1}z_{3}z_{2}z_{5}^{-1},$ \\
	&
$i_{-}(\gamma_{1})z_{5}z_{8}z_{2}z_{9}^{-1}z_{5}^{-1},\, 
 i_{-}(\gamma_{2})z_{5}z_{6}^{-1}z_{5}^{-1},\,
 i_{-}(\gamma_{3})z_{9}^{-1}z_{6}^{-1}z_{5}^{-1},\,
 i_{-}(\gamma_{4})z_{2}^{-1}z_{3}z_{1}z_{2}^{2},$ \\
	&
$i_{+}(\gamma_{1})z_{5}z_{2}z_{9}^{-1}z_{5}^{-1}, \,
 i_{+}(\gamma_{2})z_{5}z_{9}z_{6}^{-1},\,
 i_{+}(\gamma_{3})z_{2}^{-1}z_{6}^{-1},\,
 i_{+}(\gamma_{4})z_{2}^{-1}z_{1}z_{2}^{2}\,$
\\
\multicolumn{1}{c}{Torsion $\tau^{+}_{\rho_2}$} & 
$\displaystyle{-\frac{\gamma_3^2\gamma_4^5}{\gamma_2^5}
+\frac{\gamma_3^2\gamma_4^5}{\gamma_1\gamma_2^5}
+\frac{\gamma_3^2\gamma_4^6}{\gamma_2^5}}$ \\
\end{tabular}
\begin{figure*}[htbp]
\includegraphics[width=0.99\textwidth]{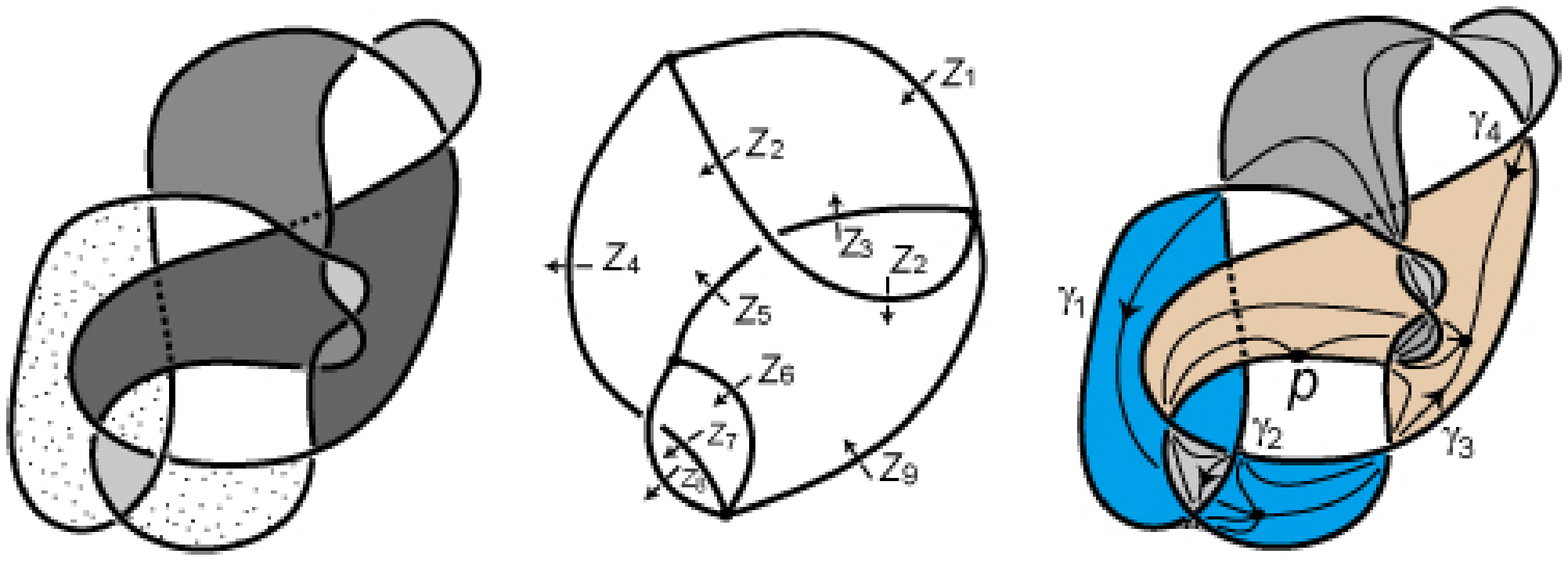}
\caption{0279}
\label{fig:0279SG}
\end{figure*}
\end{center}

\begin{center}
\begin{tabular}{lcr}
\multicolumn{2}{c}{0382}\\
\multicolumn{1}{c}{Generators} & 
$i_{-}(\gamma_{1}),\ldots, i_{-}(\gamma_{4}),\, 
 z_{1},\ldots z_{4},\,
 i_{+}(\gamma_{1}),\ldots,i_{+}(\gamma_{4})$ \\ 
\multicolumn{1}{c}{Relations} & 
$
i_{-}(\gamma_{1})z_{2}z_{1}^{-1}z_{3}z_{2}^{-1},\, 
i_{-}(\gamma_{2})z_{2}z_{3}^{-1}z_{2}z_{1}^{-2}z_{4}z_{2}^{-1},\,
i_{-}(\gamma_{3})z_{4}^{-1}z_{1}^{-1}z_{4}z_{2}^{-1},\,$ \\
	&
$
i_{-}(\gamma_{4})z_{2}^{2}z_{1}^{-1}z_{4},\,
$ \\
	&
$
i_{+}(\gamma_{1})z_{3}z_{2}^{-1}, \,
i_{+}(\gamma_{2})z_{2}z_{1}^{-2}z_{4}z_{1}^{-1},\,
i_{+}(\gamma_{3})z_{1}^{-1},\,
i_{+}(\gamma_{4})z_{4}z_{2}z_{1}^{-1}z_{4}\,
$
\\
\multicolumn{1}{c}{Torsion $\tau^{+}_{\rho_2}$} & 
$\displaystyle{\frac{1}{\gamma_1\gamma_2\gamma_4}
+\frac{1}{\gamma_1\gamma_3^2\gamma_4}
-\frac{1}{\gamma_1\gamma_3\gamma_4}}$ \\
\end{tabular}
\begin{figure*}[h]
\includegraphics[width=0.99\textwidth]{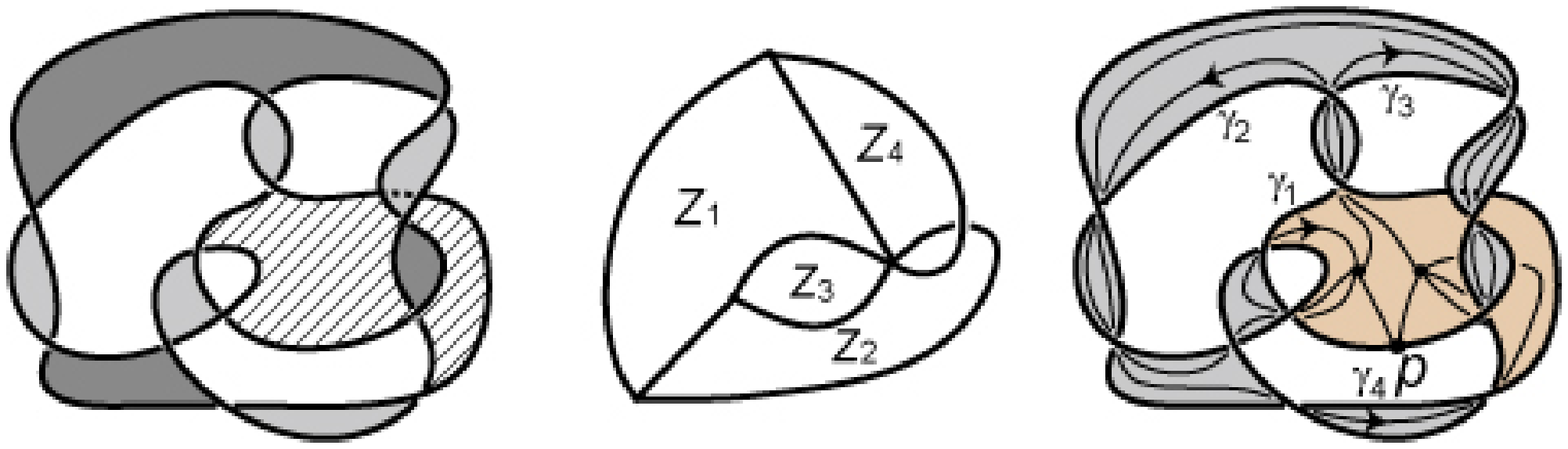}
\caption{0382}
\label{fig:0382SG}
\end{figure*}
\end{center}

\begin{center}
\begin{tabular}{lcr}
\multicolumn{2}{c}{0394}\\
\multicolumn{1}{c}{Generators} & 
$i_{-}(\gamma_{1}),\ldots, i_{-}(\gamma_{4}),\, 
 z_{1},\ldots z_{4},\,
 i_{+}(\gamma_{1}),\ldots,i_{+}(\gamma_{4})$ \\ 
\multicolumn{1}{c}{Relations} & 
$
i_{-}(\gamma_{1})z_{1}^{-1}z_{2}^{-1}z_{3},\, 
i_{-}(\gamma_{2})z_{3}^{-1}z_{4}z_{2}z_{3}z_{2}^{-1}z_{1},\,
i_{-}(\gamma_{3})z_{4}z_{2}z_{3}z_{2}^{-1}z_{1},\,
i_{-}(\gamma_{4})z_{4},\,
$ \\
	&
$
i_{+}(\gamma_{1})z_{2}^{-1}z_{3}, \,
i_{+}(\gamma_{2})z_{3}^{-1}z_{1}z_{3}^{-1}z_{4}z_{2}z_{3}z_{2}^{-1},\,
i_{+}(\gamma_{3})z_{2}z_{3}z_{2}^{-1},\,
i_{+}(\gamma_{4})z_{2}z_{4}\,
$
\\
\multicolumn{1}{c}{Torsion $\tau^{+}_{\rho_2}$} & 
$\displaystyle{\frac{1}{\gamma_1\gamma_2\gamma_3^2\gamma_4}
+\frac{1}{\gamma_1^2\gamma_2\gamma_3\gamma_4}
-\frac{1}{\gamma_1\gamma_2\gamma_3\gamma_4}}$ \\
\end{tabular}
\begin{figure*}[h]
\includegraphics[width=0.99\textwidth]{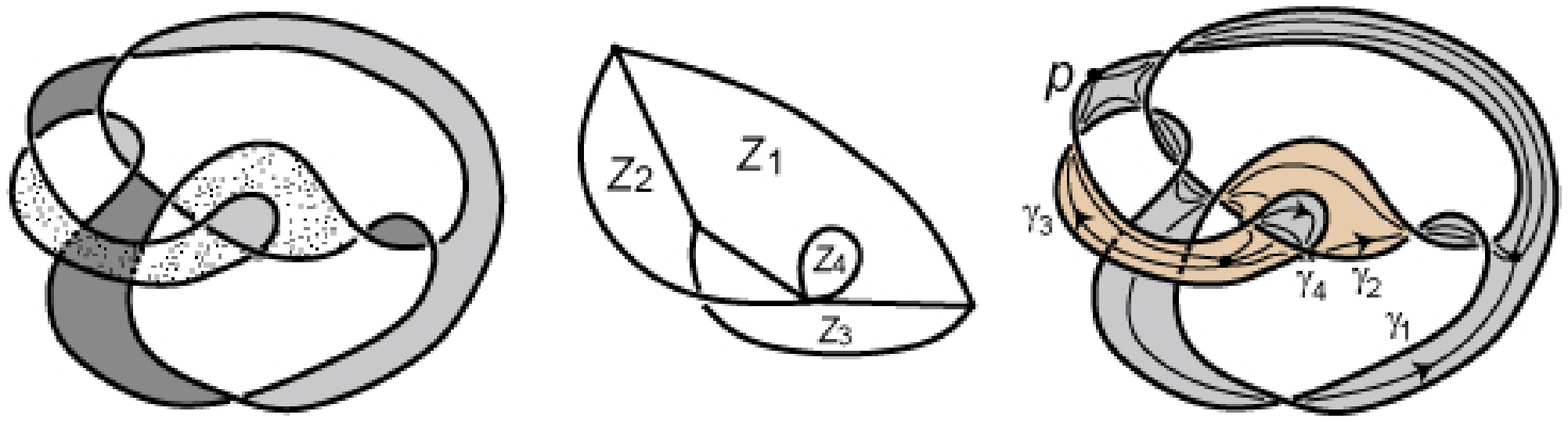}
\caption{0394}
\label{fig:0394SG}
\end{figure*}
\end{center}

\begin{center}
\begin{tabular}{lcr}
\multicolumn{2}{c}{0464}\\
\multicolumn{1}{c}{Generators} & 
$i_{-}(\gamma_1),\ldots, i_{-}(\gamma_4),\, 
z_{1},\ldots ,z_{10},\,
i_{+}(\gamma_{1}),\ldots,i_{+}(\gamma_{4})$ \\ 
\multicolumn{1}{c}{Relations} & 
$
z_{1}z_{2}z_{6}z_{7},\,z_{2}z_{9}z_{7},\,z_{3}z_{4}z_{5}z_{10}^{-1},\,
z_{4}z_{5}z_{8},\,z_{1}z_{2}z_{3}^{-1}z_{2}^{-1},\,
z_{8}z_{6}z_{8}^{-1}z_{9}^{-1}, $ \\
	&
$
i_{-}(\gamma_{1})z_{2}z_{10}z_{5}^{-1}z_{9}^{-1}z_{2}^{-1},\, 
i_{-}(\gamma_{2})z_{2}z_{10}z_{5}^{-1}z_{3}^{-1}z_{2}^{-1},\,
i_{-}(\gamma_{3})z_{2}z_{8}^{-1}z_{2}^{-1},\, 
i_{-}(\gamma_{4})z_{2}z_{1},$ \\
	&
$i_{+}(\gamma_{1})z_{2}z_{9}^{-1}z_{2}^{-1}, \,
 i_{+}(\gamma_{2})z_{2}z_{5}^{-1}z_{2}^{-1},\,
 i_{+}(\gamma_{3})z_{1}^{-1}z_{8}^{-1}z_{9}^{-1}z_{2}^{-1}z_{1},\,
 i_{+}(\gamma_{4})z_{1}^{-1}z_{7}^{-1}z_{1}\,$
\\
\multicolumn{1}{c}{Torsion $\tau^{+}_{\rho_2}$} & 
$\displaystyle{-\frac{\gamma_1^3\gamma_4^3}{\gamma_3}-\gamma_1^2\gamma_4^4
+\gamma_1^3\gamma_4^4}$ \\
\end{tabular}
\begin{figure*}[h]
\includegraphics[width=0.99\textwidth]{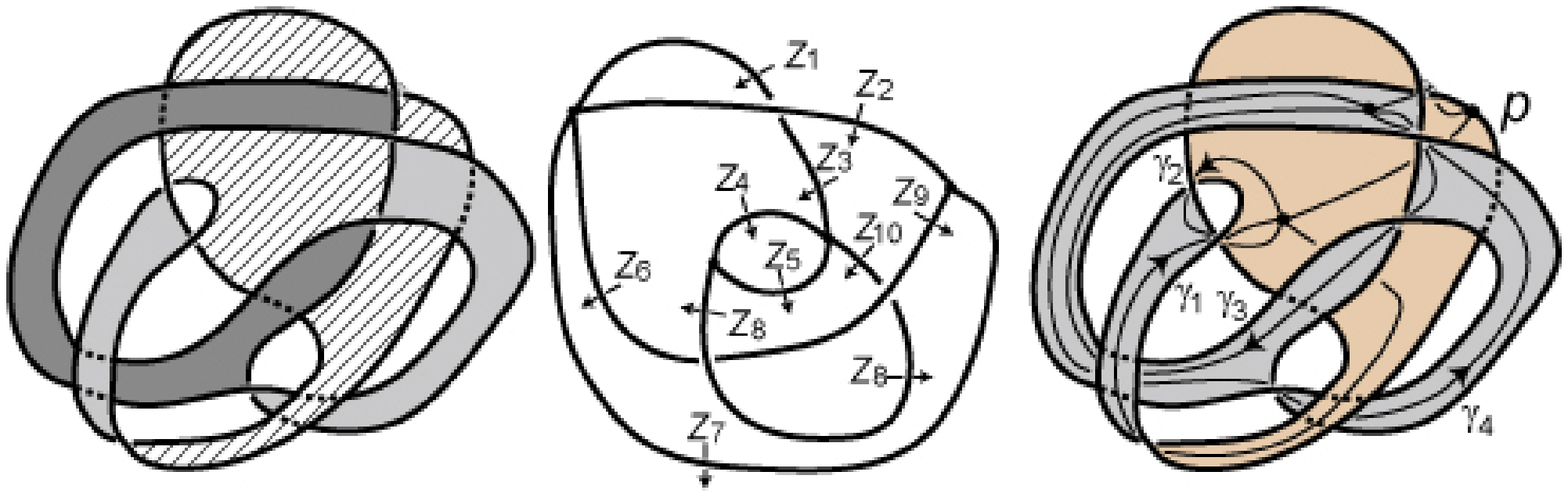}
\caption{0464}
\label{fig:0464SG}
\end{figure*}
\end{center}

\begin{center}
\begin{tabular}{lcr}
\multicolumn{2}{c}{0483}\\
\multicolumn{1}{c}{Generators} & 
$i_{-}(\gamma_1),\ldots, i_{-}(\gamma_4),\, 
z_{1},\ldots ,z_{9},\,
i_{+}(\gamma_{1}),\ldots,i_{+}(\gamma_{4})$ \\ 
\multicolumn{1}{c}{Relations} & 
$
z_{8}^{-1}z_{1}z_{4}z_{9}z_{4}^{-1},\,
z_{5}z_{6}z_{7}^{-1}z_{6}^{-1}z_{8},\,
z_{2}z_{3}z_{2}^{-1}z_{1},\,
z_{3}^{-1}z_{2}z_{3}z_{5}^{-1},\,
z_{4}z_{9}^{-1}z_{4}^{-1}z_{3},\,$ \\
	&
$
i_{-}(\gamma_{1})z_{1}z_{2}^{-1}z_{1}^{-1},\, 
i_{-}(\gamma_{2})z_{1}z_{4}^{-1}z_{8}^{-1},\,
i_{-}(\gamma_{3})z_{6}^{-1},\, 
i_{-}(\gamma_{4})z_{6}^{-1}z_{3},$ \\
	&
$i_{+}(\gamma_{1})z_{4}^{-1}z_{2}^{-1}, \,
 i_{+}(\gamma_{2})z_{4}^{-1},\,
 i_{+}(\gamma_{3})z_{5}z_{6}^{-1}z_{8},\,
 i_{+}(\gamma_{4})z_{8}^{-1}z_{3}$
\\
\multicolumn{1}{c}{Torsion $\tau^{+}_{\rho_2}$} & 
$\displaystyle{\frac{1}{\gamma_1\gamma_3\gamma_4^2}
-\frac{\gamma_2}{\gamma_1^2\gamma_3\gamma_4^2}
-\frac{1}{\gamma_1\gamma_3\gamma_4}}$ \\
\end{tabular}
\begin{figure*}[h]
\includegraphics[width=0.99\textwidth]{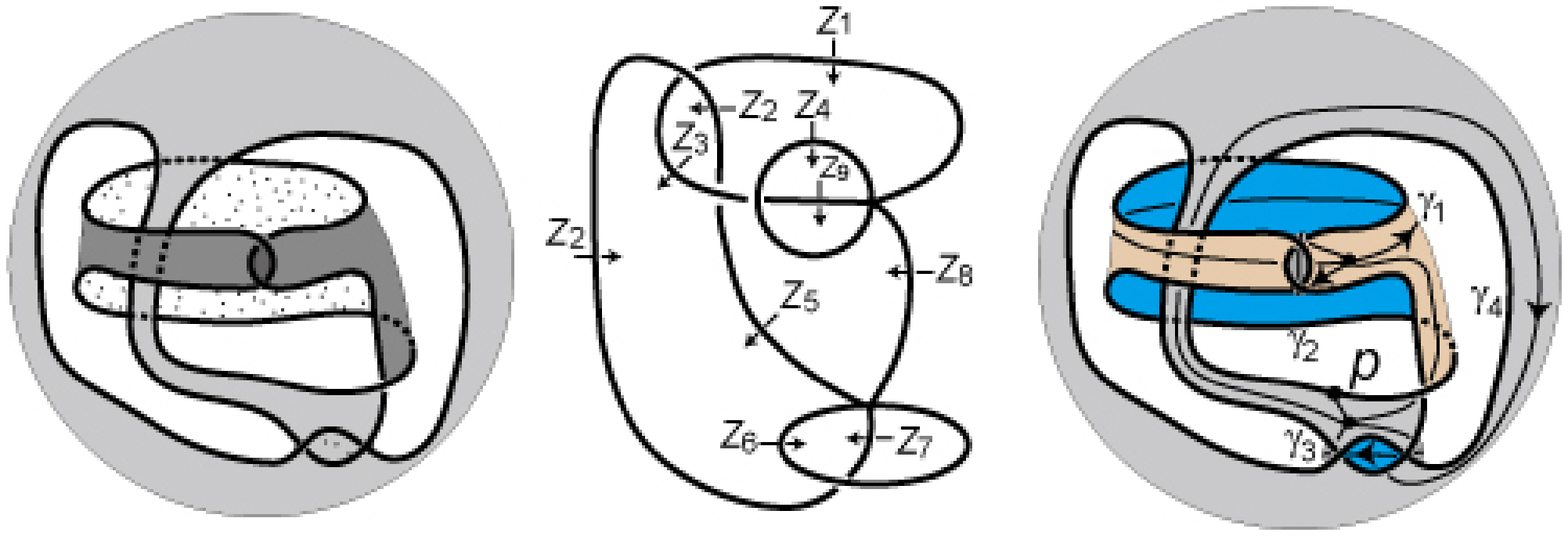}
\caption{0483}
\label{fig:0483SG}
\end{figure*}
\end{center}

\begin{center}
\begin{tabular}{lcr}
\multicolumn{2}{c}{0535}\\
\multicolumn{1}{c}{Generators} & 
$i_{-}(\gamma_1),\ldots, i_{-}(\gamma_4),\, 
z_{1},\ldots ,z_{10},\,
i_{+}(\gamma_{1}),\ldots,i_{+}(\gamma_{4})$ \\ 
\multicolumn{1}{c}{Relations} & 
$
z_{1}z_{2}z_{3},\,z_{2}z_{6}z_{7}^{-1},\,
z_{10}^{-1}z_{4}z_{5}z_{8}^{-1}z_{7},\,
z_{1}z_{10}z_{9},\,z_{2}z_{3}z_{2}^{-1}z_{4}^{-1},\,
z_{2}z_{6}^{-1}z_{2}^{-1}z_{5}^{-1},$ \\
	&
$
i_{-}(\gamma_{1})z_{10}^{-1},\, 
i_{-}(\gamma_{2})z_{10}^{-1}z_{1}^{-1}z_{3}^{-1}z_{1}^{-1}z_{10},\,
i_{-}(\gamma_{3})z_{7}^{-1}z_{1}^{-1}z_{10},\, 
i_{-}(\gamma_{4})z_{6}z_{3}^{-1}z_{7},$ \\
	&
$i_{+}(\gamma_{1})z_{7}^{-1}z_{9}, \,
 i_{+}(\gamma_{2})z_{7}^{-1}z_{1}^{-1}z_{3}^{-1}z_{10}z_{7},\,
 i_{+}(\gamma_{3})z_{7}^{-1}z_{3}z_{10}z_{7},\,
 i_{+}(\gamma_{4})z_{7}^{-1}z_{6}z_{3}^{-1}z_{7}$
\\
\multicolumn{1}{c}{Torsion $\tau^{+}_{\rho_2}$} & 
$\displaystyle{-\frac{1}{\gamma_1^{11}\gamma_2^6\gamma_3^6\gamma_4^{15}}
+\frac{1}{\gamma_1^{10}\gamma_2^5\gamma_3^6\gamma_4^{15}}
-\frac{1}{\gamma_1^{10}\gamma_2^5\gamma_3^6\gamma_4^{14}}}$ \\
\end{tabular}
\begin{figure*}[h]
\includegraphics[width=0.99\textwidth]{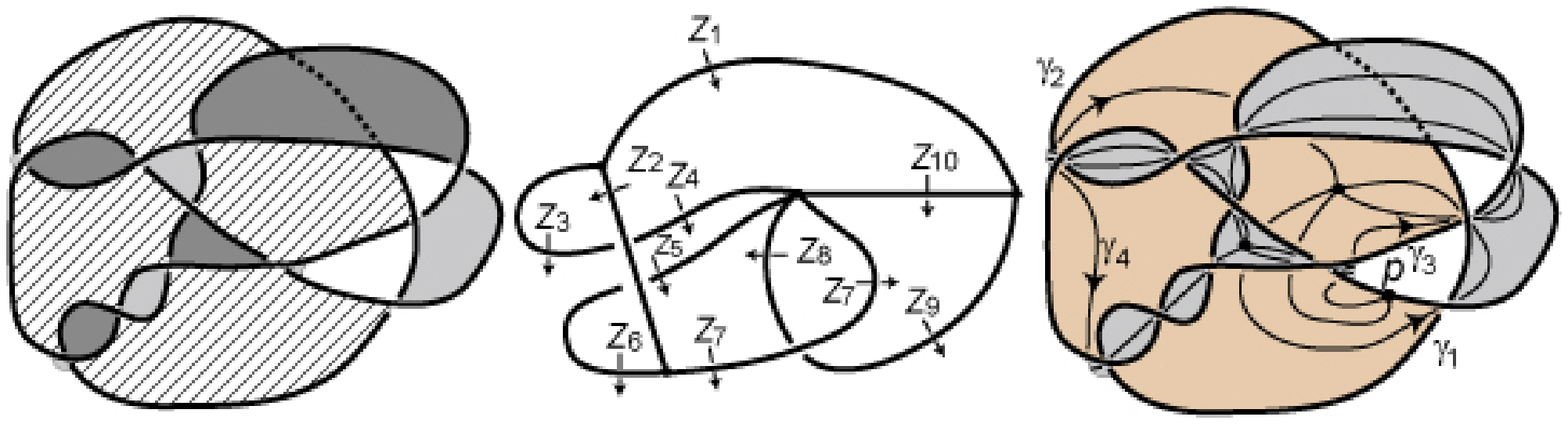}
\caption{0535}
\label{fig:0535SG}
\end{figure*}
\end{center}

\begin{center}
\begin{tabular}{lcr}
\multicolumn{2}{c}{0650}\\
\multicolumn{1}{c}{Generators} & 
$i_{-}(\gamma_1),\ldots, i_{-}(\gamma_4),\, 
z_{1},\ldots ,z_{11},\,
i_{+}(\gamma_{1}),\ldots,i_{+}(\gamma_{4})$ \\ 
\multicolumn{1}{c}{Relations} & 
$
z_{2}z_{3}z_{1}^{-1}z_{4}z_{1},\,
z_{2}z_{6}z_{8}^{-1}z_{6}^{-1}z_{11}^{-1},\,
z_{1}z_{5}^{-1}z_{1}^{-1}z_{4},\,
z_{3}z_{6}z_{9}^{-1}z_{6}^{-1},\,
z_{9}z_{8}^{-1}z_{7}z_{8},\,$ \\
	&
$
z_{8}z_{7}z_{10}^{-1}z_{7}^{-1},\,
z_{10}z_{6}^{-1}z_{11}z_{6}, $ \\
	&
$
i_{-}(\gamma_{1})z_{2}z_{6}^{-1}z_{2}^{-1},\, 
i_{-}(\gamma_{2})z_{2}z_{7}z_{6}^{-1}z_{2}^{-1},\,
i_{-}(\gamma_{3})z_{6}z_{8}z_{6}^{-1}z_{2}^{-1},\, 
i_{-}(\gamma_{4})z_{2}z_{3}z_{1}^{-1},$ \\
	&
$i_{+}(\gamma_{1})z_{11}z_{6}^{-1}z_{2}^{-1}, \,
 i_{+}(\gamma_{2})z_{2}z_{3}^{-1}z_{2}^{-1},\,
 i_{+}(\gamma_{3})z_{1}z_{6}z_{8}z_{6}^{-1},\,
 i_{+}(\gamma_{4})z_{1}^{-1}$
\\
\multicolumn{1}{c}{Torsion $\tau^{+}_{\rho_2}$} &
$\displaystyle{\frac{1}{\gamma_1\gamma_2^3\gamma_3^2\gamma_4^2}
-\frac{1}{\gamma_1\gamma_2^3\gamma_3\gamma_4}
+\frac{1}{\gamma_1\gamma_2^2\gamma_3\gamma_4}}$ \\
\end{tabular}
\begin{figure*}[h]
\includegraphics[width=0.99\textwidth]{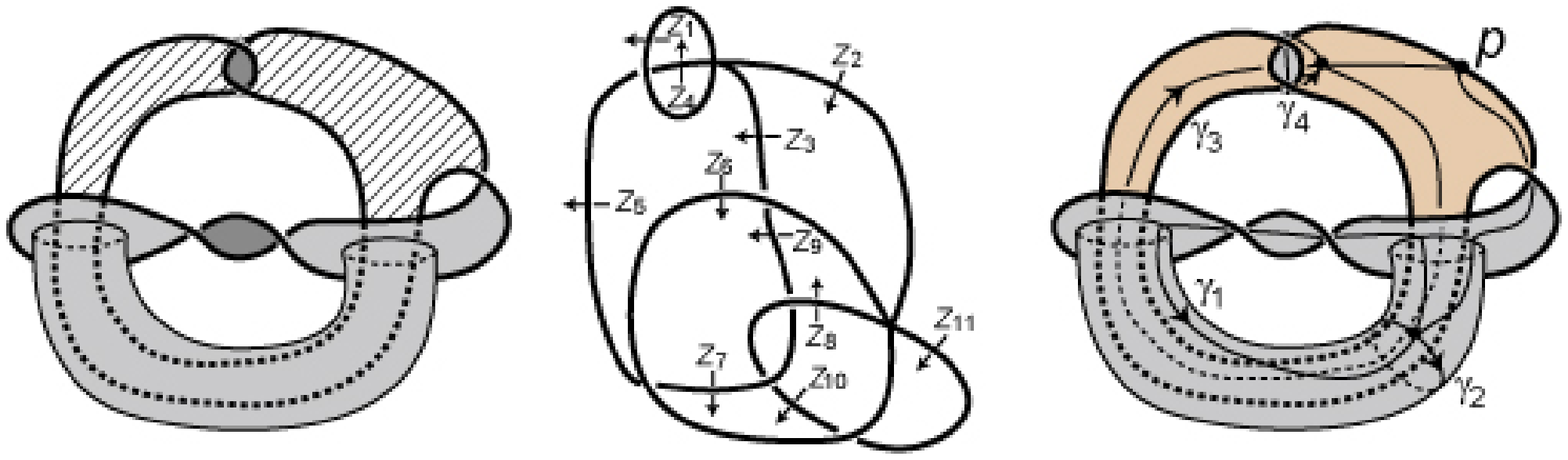}
\caption{0650}
\label{fig:0650SG}
\end{figure*}
\end{center}

\begin{center}
\begin{tabular}{lcr}
\multicolumn{2}{c}{0801}\\
\multicolumn{1}{c}{Generators} & 
$i_{-}(\gamma_1),\ldots, i_{-}(\gamma_4),\, 
z_{1},\ldots ,z_{9},\,
i_{+}(\gamma_{1}),\ldots,i_{+}(\gamma_{4})$ \\ 
\multicolumn{1}{c}{Relations} & 
$
z_{1}^{-1}z_{6}z_{7}z_{8}^{-1}z_{9}^{-1},\,
z_{3}z_{4}z_{9}z_{6}^{-1},\,
z_{2}z_{4}z_{5},\,
z_{2}z_{6}z_{7}^{-1}z_{6}^{-1},\,
z_{2}z_{3}^{-1}z_{2}^{-1}z_{1}, $ \\
	&
$
i_{-}(\gamma_{1})z_{6}z_{7}z_{8}^{-1}z_{6},\, 
i_{-}(\gamma_{2})z_{1}z_{2}z_{8}z_{7}^{-1}z_{6}^{-1},\,
i_{-}(\gamma_{3})z_{9}z_{6}^{-1}z_{2}^{-1},\, 
i_{-}(\gamma_{4})z_{5}^{-1}z_{9}^{-1}z_{5}^{-1},$ \\
	&
$i_{+}(\gamma_{1})z_{6}z_{9}, \,
 i_{+}(\gamma_{2})z_{6}z_{2}z_{6}z_{9}^{-1}z_{6}^{-1},\,
 i_{+}(\gamma_{3})z_{5}z_{9}z_{6}^{-1},\,
 i_{+}(\gamma_{4})z_{4}z_{9}^{-1}z_{5}^{-1}$
\\
\multicolumn{1}{c}{Torsion $\tau^{+}_{\rho_2}$} & 
$-\gamma_1^2\gamma_3^2\gamma_4
+\gamma_1^2\gamma_2\gamma_3^2\gamma_4
-\gamma_1^2\gamma_2\gamma_3^3\gamma_4^2$ \\
\end{tabular}
\begin{figure*}[h]
\includegraphics[width=0.99\textwidth]{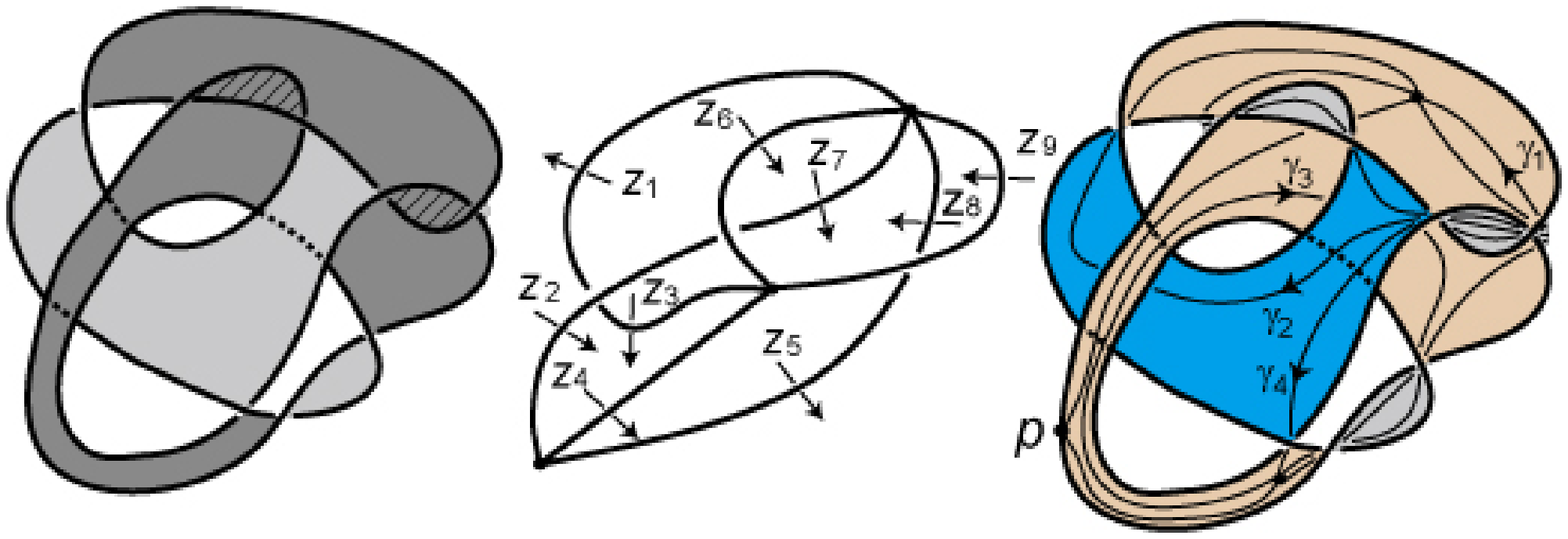}
\caption{0801}
\label{fig:0801SG}
\end{figure*}
\end{center}

\begin{center}
\begin{tabular}{lcr}
\multicolumn{2}{c}{0815}\\
\multicolumn{1}{c}{Generators} & 
$i_{-}(\gamma_1),\ldots, i_{-}(\gamma_4),\, 
z_{1},\ldots ,z_{11},\,
i_{+}(\gamma_{1}),\ldots,i_{+}(\gamma_{4})$ \\ 
\multicolumn{1}{c}{Relations} & 
$
z_{1}z_{9}z_{6},\,
z_{1}z_{2}^{-1}z_{4}^{-1},\,
z_{4}z_{11}^{-1}z_{5},\,
z_{10}^{-1}z_{5}^{-1}z_{6}z_{7}z_{8},\,
z_{8}^{-1}z_{6}^{-1}z_{9}z_{6},\,$ \\
	&
$
z_{7}^{-1}z_{6}^{-1}z_{3}z_{6},\,
z_{4}z_{3}^{-1}z_{4}^{-1}z_{10}, $ \\
	&
$
i_{-}(\gamma_{1})z_{4}z_{3}^{-1}z_{4}^{-1},\, 
i_{-}(\gamma_{2})z_{4}z_{11},\,
i_{-}(\gamma_{3})z_{9},\, 
i_{-}(\gamma_{4})z_{2}^{-1}z_{9}^{-1},$ \\
	&
$i_{+}(\gamma_{1})z_{2}^{-1}z_{3}^{-1}z_{4}^{-1}, \,
 i_{+}(\gamma_{2})z_{11}z_{1},\,
 i_{+}(\gamma_{3})z_{9}z_{3}^{-1}z_{1},\,
 i_{+}(\gamma_{4})z_{9}z_{2}^{-1}z_{9}^{-1}$
\\
\multicolumn{1}{c}{Torsion $\tau^{+}_{\rho_2}$} & 
$\displaystyle{-\frac{\gamma_1^3\gamma_2^5}{\gamma_4^6}
+\frac{\gamma_1^2\gamma_2^4}{\gamma_4^5}
+\frac{\gamma_1^3\gamma_2^5}{\gamma_4^5}}$ \\
\end{tabular}
\begin{figure*}[h]
\includegraphics[width=0.99\textwidth]{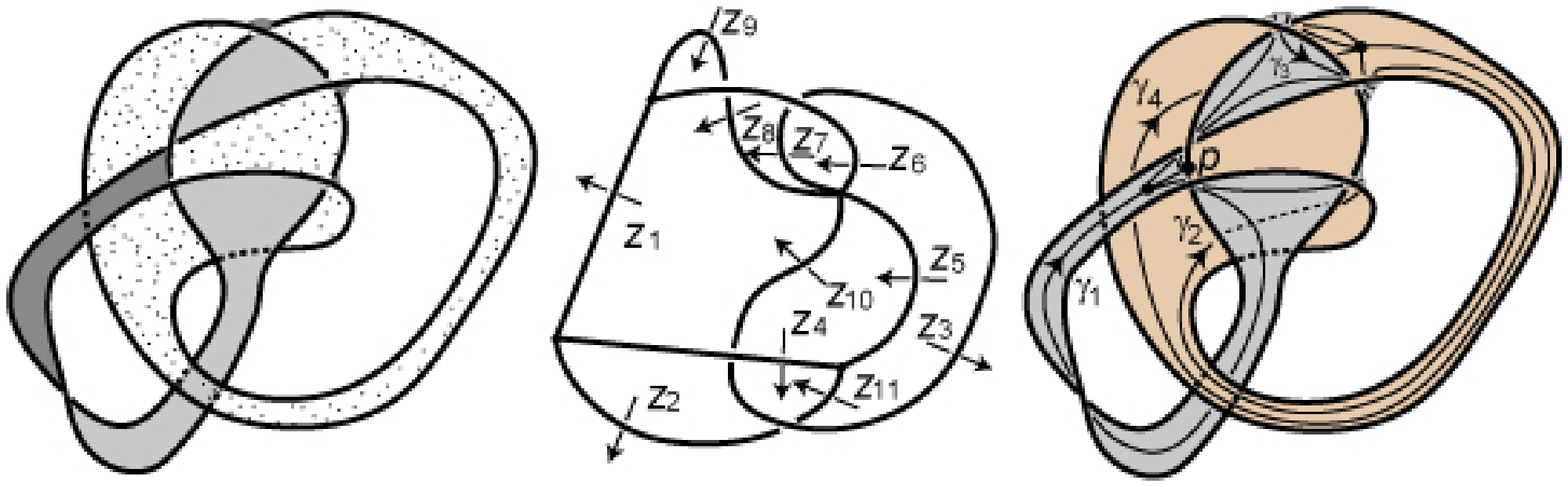}
\caption{0815}
\label{fig:0815SG}
\end{figure*}
\end{center}

\begin{remark}
According to \cite{bg} and \cite{juhasz2}, 
these knots have unique minimal genus Seifert surfaces.
\end{remark}


\section{Magnus matrix and concordances of Seifert surfaces}\label{sec:magnus}

Not only the torsion $\tau_\rho^+$ but 
the Magnus matrix $r_\rho$ can be used 
as a fibering obstruction of homologically fibered knots. 
In fact, for a fibered knot $K$ with its unique minimal genus 
Seifert surface $R$ of genus $g$, 
the sutured manifold $M_R$ is given by 
a mapping cylinder 
$(\Sigma_{g,1} \times [0,1], \mathrm{id} \times 1, 
\varphi \times 0)$ of $\varphi$, which is an element 
of the mapping class group of $\Sigma_{g,1}$ and is 
uniquely (not up to conjugation) 
determined after fixing 
an identification of $R$ with $\Sigma_{g,1}$. 
Then the Magnus matrix $r_\rho (M_R)$ 
associated with a homomorphism 
$\rho:G(K) \to \Gamma$ is given by 
$\overline{
\sideset{^{\rho \circ i_+}\!}{}
{\mathop{\left({\displaystyle\frac{\partial 
\varphi(\gamma_j)}{\partial \gamma_i}} \right)}\nolimits}
}_{1 \le i,j \le 2g}$. 
In particular, 
all the entries are elements of $\mathbb{Z} \Gamma$. 
Therefore, for the detection of non-fiberedness of 
a non-fibered homologically fibered knot $K$, 
it suffices to find 
a minimal genus Seifert surface $R$, 
an identification of $R$ with $\Sigma_{g,1}$ and 
a homomorphism $\rho:G(K) \to \Gamma$ to a PTFA group $\Gamma$ 
whose Magnus matrix has an entry not contained in 
$\mathbb{Z} \Gamma$.

\begin{example}
If we continue the computation for the knot 0057 
in Section \ref{sec:sample}, 
we can see that the $(1,3)$-entry of $r_{\rho_2} (M_R)$ is 
$\displaystyle\frac{\gamma_4}{1+\gamma_2-\gamma_2 \gamma_4}$, 
not an element of $\mathbb{Z} D_2 (K)$. 
This shows that the knot 0057 is not fibered. 
\end{example}

In the usage of $r_\rho$ as a fibering obstruction, 
its invariance under homology cobordisms of 
homology cylinders is convenient. We first recall the definition of 
homology cobordisms of homology cylinders: 
\begin{definition}
Two homology cylinders $M=(M,i_+,i_-), 
N=(N,j_+,j_-) \in \mathcal{C}_{g,1}$ 
are said to be {\it $($rational\/$)$ 
homology cobordant} if there exists 
a smooth compact oriented 4-manifold $W$ such that: 
\begin{enumerate}
\item $\partial W = M \cup (-N) /(i_+ (x)= j_+(x) , \,
i_- (x)=j_-(x)) \quad x \in \Sigma_{g,1}$; and 
\item the inclusions $M \hookrightarrow W$, $N \hookrightarrow W$ 
induce isomorphisms on the (rational) homology,  
\end{enumerate}
where $-N$ denotes $N$ with opposite orientation. 
\end{definition}
\begin{proposition}
\label{prop:h-inv}
Suppose that $M=(M,i_+,i_-), N=(N,j_+,j_-) \in 
\mathcal{C}_{g,1}^\mathbb{Q}$ 
are rational homology cobordant by a rational homology cobordism $W$. 
Let $\rho: \pi_1 (W) \to \Gamma$ be a homomorphism to a PTFA group. 
Then the Magnus matrices $r_\rho (M)$ and $r_\rho (N)$ 
associated with 
$\pi_1 (M) \to \pi_1 (W) \xrightarrow{\rho} \Gamma$ 
and $\pi_1 (N) \to \pi_1 (W) \xrightarrow{\rho} \Gamma$ 
are defined, and 
$r_\rho (M) = r_\rho (N)$ holds. 
\end{proposition}
\begin{proof}
We can apply the argument of \cite[Section 3.1]{sakasai08} and 
we omit the details. 
\end{proof}
To interpret the homology cobordant relation in terms of homologically 
fibered knots, we introduce {\it concordances of Seifert surfaces} 
defined by Myers. 
\begin{definition}[Myers \cite{myers}]
Seifert surfaces $R$, $R'$ of genus $g$ 
for knots $K$, $K'$ in $S^3$ are said to be 
{\it concordant} if there is a smooth embedding 
$I: \Sigma_{g,1} \times [0,1] \to S^3 \times [0,1]$ such that 
$I (\Sigma_{g,1} \times \{ 0 \}) = R$ and 
$I (\Sigma_{g,1} \times \{ 1 \}) = R'$. 
\end{definition}
Using this terminology, we have the following relationship between 
concordances of Seifert surfaces for homologically fibered knots 
and homology cylinders. 
\begin{proposition}\label{prop:concordant}
Let $K$ be a $($rationally\/$)$ homologically fibered knot 
with a minimal genus Seifert surface $R$ of genus $g$. 
Suppose $R$ is concordant to 
another Seifert surface $R'$ of a knot $K'$. Then 
$K'$ is also a $($rationally\/$)$ homologically fibered knot 
of genus $g$ with a minimal genus Seifert surface $R'$ such that 
$M_R$ and $M_{R'}$ are $($rational\/$)$ homology cobordant 
as $($rational\/$)$ homology cylinders.
\end{proposition}
\begin{proof}
Let $W$ be a manifold obtained from $S^3 \times [0,1]$ by cutting 
along the image of 
an embedding $I: \Sigma_{g,1} \times [0,1] \to S^3 \times [0,1]$ 
which connects $R$ and $R'$. 
Then it is straightforward to check our assertions by 
observing the Mayer-Vietoris exact sequence of 
$S^3 \times [0,1]=W \cup I(\Sigma_{g,1} \times [0,1])$ with 
the intersection homeomorphic to 
$(\partial M_R) \times [0,1] = (\Sigma_{g,1} \cup (-\Sigma_{g,1})) 
\times [0,1]$. We omit the details.  
\end{proof}
The following theorem enables us to produce 
infinitely many Seifert surfaces which are 
concordant to a given one. 
\begin{theorem}[Myers \cite{myers}]\label{thm:myers}
If a Seifert surface $R$ of a knot $K$ is not a disk, then $R$ is 
concordant to $R'$ such that: 
\begin{itemize}
\item[$(1)$] $K'=\partial R'$ is hyperbolic; and
\item[$(2)$] $S^3 - K'$ has arbitrarily large hyperbolic volume. 
\end{itemize}
\end{theorem}

For $(M,i_+,i_-) \in \mathcal{C}_{g,1}$, consider 
a homomorphism \[\rho_M: \pi_1 (M) \to H_1 (M) \xrightarrow{i_+^{-1}} 
H_1 (\Sigma_{g,1}).\]
We can easily check that if $W$ is a homology cobordism between 
$M$ and $N$ in $\mathcal{C}_{g,1}$, then 
there exists an extension 
$\rho_W:\pi_1 (W) \to H_1 (\Sigma_{g,1})$ 
of $\rho_M$ and $\rho_N$. 
Note that $\rho_M$ can be regarded as 
a restriction of $\rho_2$ when $M$ is obtained from 
a homologically fibered knot (recall the exact sequence 
(\ref{eq:seq})). Consequently 
we can combine Theorem \ref{thm:myers} with 
Proposition \ref{prop:h-inv} as follows: 
\begin{theorem}\label{thm:detect}
Let $K$ be a homologically fibered knot 
with a minimal genus Seifert surface $R$. 
If $K$ is shown to be non-fibered by using 
$r_{\rho_2} (M_R) (=r_{\rho_{M_R}} (M_R))$, 
then $K'=\partial R'$ is also non-fibered for any 
Seifert surface $R'$ concordant to $R$. Moreover, there exist 
infinitely many such $K'= \partial R'$. 
\end{theorem}
We may take $K$ to be a homologically fibered knot in Section 
\ref{sec:sample}. 
Then Theorem \ref{thm:detect} shows that there does exist infinitely 
many homologically fibered knots whose non-fiberedness are detected by 
the Magnus matrices. 

\begin{example}\label{ex:concordance}
Let $K$ be the knot as the boundary of the Seifert surface $R$ 
illustrated in Figure \ref{fig:concordance}. 

\begin{figure}[htbp]
\begin{center}
\includegraphics[width=.8\textwidth]{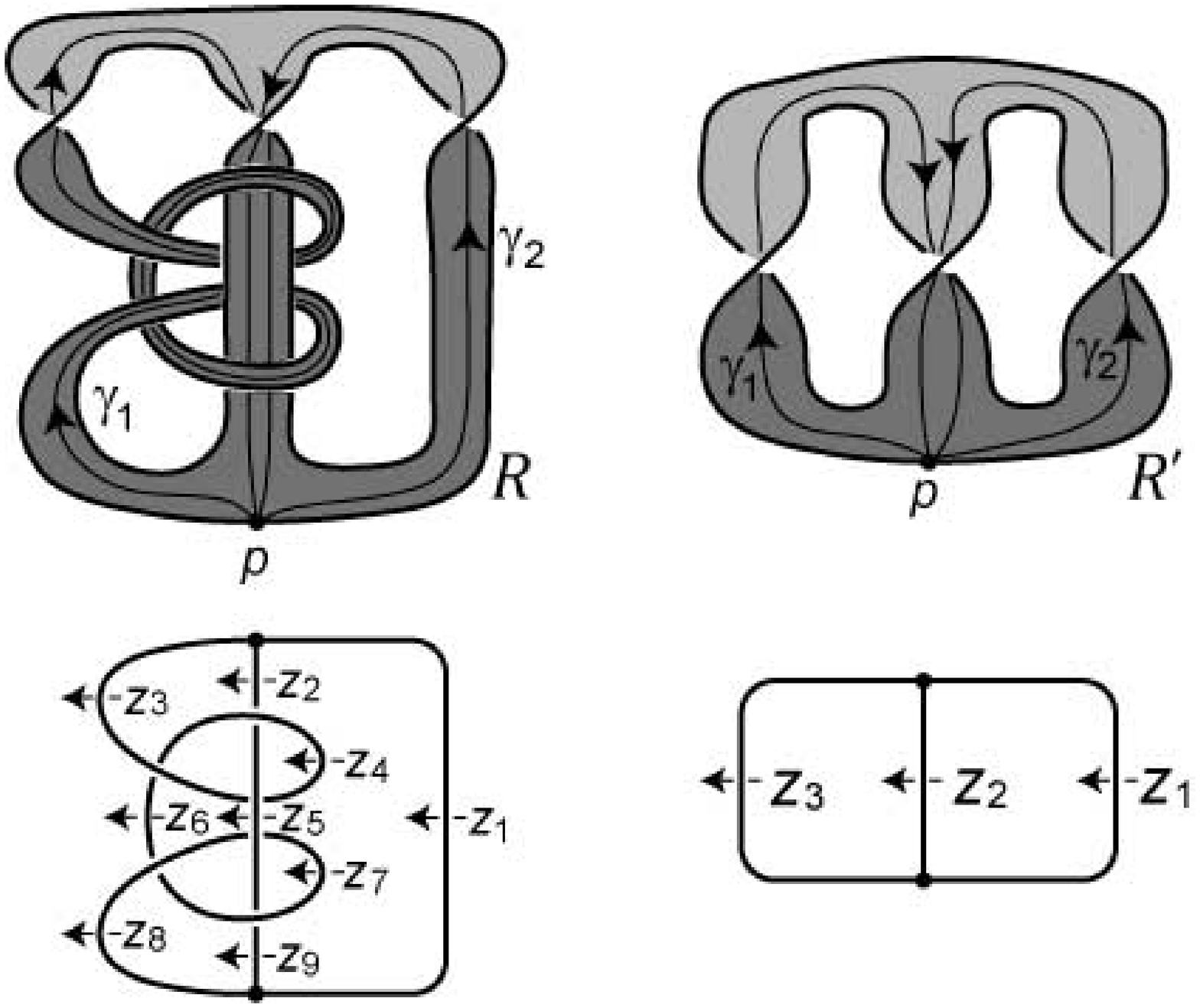}
\end{center}
\caption{Concordant Seifert surfaces}
\label{fig:concordance}
\end{figure}

$R$ is concordant to the minimal genus Seifert surface $R'$ 
of the trefoil knot, which is fibered. 
Proposition \ref{prop:concordant} shows that $K$ is 
a homologically fibered knot and $R$ is of minimal genus. 
An admissible presentation of $\pi_1 (M_R)$ is given by
{\small \begin{center}
\begin{tabular}{lcr}
\multicolumn{1}{c}{Generators} & 
$i_{-}(\gamma_1),\,i_{-}(\gamma_2),\, 
z_{1},\ldots ,z_{9},\,
i_{+}(\gamma_{1}),\,i_{+}(\gamma_2)$ \\ 
\multicolumn{1}{c}{Relations} & 
$z_{1}z_{2}z_{3},\,z_{1}z_{9}z_{8},\,
z_{4}z_{5}z_{4}^{-1}z_2^{-1},\,
z_{4}^{-1}z_{5}z_{3}^{-1}z_5^{-1},\,
z_{3}z_{6}z_{3}^{-1}z_4,\,
z_{7}z_{5}z_{8}z_{5}^{-1}$,\\
& $z_{7}^{-1}z_{9}z_{7}z_{5}^{-1},\, 
i_{-}(\gamma_{1})z_{1}z_7z_{4}^{-1}z_2z_5^{-1}z_3z_8^{-1}z_5,\, 
i_{-}(\gamma_{2})z_{8}^{-1}z_7z_4^{-1}z_1^{-1},$ \\
	&
$i_{+}(\gamma_{1})z_7 z_{4}^{-1} z_2 z_5^{-1} z_3 z_8^{-1} z_5, \ 
i_{+}(\gamma_{2})z_{7}z_4^{-1}z_{1}^{-1}.$
\end{tabular}
\end{center}
}
\noindent
From this, we have 
\begin{align*}
\det (\tau_{\rho_2}^+ (M_R)) &= 3-\frac{1}{\gamma_1}-\gamma_1
-\frac{\gamma_1}{\gamma_2}+\frac{\gamma_1^2}{\gamma_2} 
+\frac{\gamma_2}{\gamma_1^2}-\frac{\gamma_2}{\gamma_1},\\
r_{\rho_2} (M_R) &= 
\begin{pmatrix}
1 & \gamma_2^{-1}\\
-\gamma_1^{-1} \gamma_2 & 1-\gamma_1^{-1}
\end{pmatrix}.
\end{align*}
\noindent
On the other hand, an admissible presentation of 
$\pi_1 (M_{R'})$ is given by 
{\small \begin{center}
\begin{tabular}{lcr}
\multicolumn{1}{c}{Generators} & 
$i_{-}(\gamma_1),\,i_{-}(\gamma_2),\, 
z_{1},\,z_2,\,z_3,\,
i_{+}(\gamma_{1}),\,i_{+}(\gamma_{2})$ \\ 
\multicolumn{1}{c}{Relations} & 
$z_{1}z_{2}z_{3},\,i_-(\gamma_1)z_{3}^{-1},\,
i_-(\gamma_2)z_3^{-1}z_1^{-1},\,i_+(\gamma_1) z_2,\,
i_+(\gamma_2)z_1^{-1}$ 
\end{tabular}
\end{center}
}
\noindent
and we have
\begin{align*}
\det (\tau_{\rho_2}^+ (M_R)) &= \frac{1}{\gamma_2},\\
r_{\rho_2} (M_R) &= 
\begin{pmatrix}
1 & \gamma_2^{-1}\\
-\gamma_1^{-1} \gamma_2 & 1-\gamma_1^{-1}
\end{pmatrix}.
\end{align*}
\end{example}

\begin{remark}
As seen in Example \ref{ex:concordance}, 
$\Gamma$-torsion generally changes 
under homology cobordisms. 
However, Cha-Friedl-Kim \cite{cfk} recently found 
a way to extract homology cobordant invariants from 
the torsion $\tau_{\rho_2}^+$ by taking 
a certain quotient of the target. Then they applied it to 
the {\it homology cobordism group of 
homology cylinders} and showed that this group has $\mathbb{Z}_2^\infty$ 
as an abelian quotient. By Proposition \ref{prop:unique}, 
we may regard this abelian quotient as an invariant of 
homologically fibered knots. In fact, it is unchanged under 
concordances of Seifert surfaces. 
\end{remark}


\section{MATHEMATICA program}\label{sec:program}

The following is a MATHEMATICA program which 
calculates the invariants discussed in Sections \ref{sec:sample}, 
\ref{sec:HFK12} and \ref{sec:magnus}. 

\medskip

{\scriptsize
\begin{verbatim}
h1Class = {};
h1Monodromy = {};
torsionMatrix = {};
magnusMatrix = {};

invariants[g_, z_, RELATIONS_] := 
  Module[{reindexedRel, h1Matrix, i, alex},
   GENUS = g;
   Ztotal = z;
   
   reindexedRel = Map[reindexing, RELATIONS, {2}];
   
   h1Matrix = -Map[Take[#, -2 GENUS] &, homologyComputation[reindexedRel]];
   h1Class = 
    Join[Map[monomialExpression, h1Matrix], 
     Table[ToExpression[ToString[SequenceForm["\[Gamma]", i]]], {i, 2 GENUS}]];
   Print["Homology classes of generators = ", h1Class // DisplayForm];
   
   h1Monodromy = Transpose[Take[h1Matrix, 2 GENUS]];
   Print["Homological monodromy = ", h1Monodromy // MatrixForm];
   
   alex = Transpose[makeAlexanderMatrix[reindexedRel]];
   torsionMatrix = Take[alex, 2 GENUS + Ztotal];
   Print["torsion matrix = ", torsionMatrix // MatrixForm];
   Print["det(torsion) = ", Expand[Det[torsionMatrix]]];
   
   magnusMatrix = Simplify[Transpose[
     Take[Transpose[-Drop[alex, 2 GENUS + Ztotal].Inverse[
           torsionMatrix]], 2 GENUS]]];
   Print["Magnus matrix = ", magnusMatrix // MatrixForm]
   ];

reindexing[num_] := 
  Module[{numString, sg}, 
    If[NumberQ[num], num + 2 GENUS*Sign[num],
     numString = ToString[num];
     sg = If[StringTake[numString, 1] == "-", 1, 0];
     If[StringTake[numString, {1 + sg}] == "m", 
      ((-1)^sg)*ToExpression[StringDrop[numString, 1 + sg]], 
      ((-1)^sg)*(ToExpression[StringDrop[numString, 1 + sg]] + 2 GENUS + Ztotal)]]
   ];

homologyComputation[rel_] := 
  Module[{i, j}, 
   RowReduce[Table[Count[rel[[i]], j] - Count[rel[[i]], -j], 
   {i, 1, 2 GENUS + Ztotal}, {j, 1, 4 GENUS + Ztotal}]]];

monomialExpression[list_] := 
  Module[{i, prod = 1}, 
   For[i = 1, i <= 2 GENUS, i++, 
    prod = prod*(ToExpression[ToString[SequenceForm["\[Gamma]", i]]]^list[[i]])];
   prod];

makeAlexanderMatrix[rel_] := 
  Module[{i, j}, 
   Table[foxDer[rel[[i]], j], {i, 1, Length[rel]}, {j, 1, 4 GENUS + Ztotal}]];

foxDer[word_, var_] := 
  Module[{entry = 0, i}, 
   For[i = 1, i <= Length[word], i++, 
    Which[word[[i]] == var, 
     entry = entry + (makeMonomial[Take[word, i - 1]]^(-1)), 
     word[[i]] == -var, 
     entry = entry - (makeMonomial[Take[word, i]]^(-1))]];
   entry];

makeMonomial[list_] := 
  Module[{prod = 1}, 
   For[i = 1, i <= Length[list], i++, 
    prod = prod*(h1Class[[Abs[list[[i]]]]]^Sign[list[[i]]])];
   prod];

\end{verbatim}
}

A computation by this program goes as follows. 
Let $(M,i_+,i_-) \in \mathcal{C}_{g,1}$ with an admissible presentation 
\[\langle i_- (\gamma_1),\ldots,i_- (\gamma_{2g}), 
z_1 ,\ldots, z_l, 
i_+ (\gamma_1),\ldots,i_+ (\gamma_{2g}) \mid 
r_1, \ldots, r_{2g+l}
\rangle\]
of $\pi_1 (M)$. 
The main function in the program is $\mathtt{invariants}$ having 
three slots as the input. These slots correspond to the genus $g$, 
the number $l$ of $z$-generators and the list of relations. 
For each word in the relations, we make a list by replacing 
$i_-(\gamma_j)^{\pm 1}$, $z_j^{\pm 1}$ and $i_+(\gamma_j)^{\pm 1}$ by 
$\pm \mathtt{mj}$, $\pm \mathtt{j}$ and $\pm \mathtt{pj}$. 
By lining up them, 
we obtain the list of relations. 

When we compute the case of the knot 0815 
with an admissible presentation of $\pi_1 (M_R)$ 
of the sutured manifold $M_R$ as in Section \ref{sec:HFK12}, 
for example, the input is:
{\small 
\begin{verbatim}
invariants[2, 11, {{1, 9, 6}, {1, -2,-4}, {4,-11, 5}, 
  {-10, -5, 6, 7, 8}, {-8, -6, 9, 6}, {-7, -6, 3, 6}, 
  {4, -3, -4, 10}, {m1, 4, -3, -4}, {m2, 4, 11}, 
  {m3, 9}, {m4, -2, -9}, {p1, -2, -3, -4}, {p2, 11, 1}, 
  {p3, 9, -3, 1}, {p4, 9, -2, -9}}]
\end{verbatim}}
Then the function returns homology classes of generators in terms of 
$\mathtt{\gamma j}:=i_+(\gamma_j) \in H_1 (M_R)$, 
the homological monodromy matrix $\sigma (M_R)$, 
the torsion matrix $\tau_{\rho_2}^+ (M_R)$ and 
the Magnus matrix $r_{\rho_2}(M_R)$. These data can be referred 
as the variables $\mathtt{h1Class}$, $\mathtt{h1Monodromy}$, 
$\mathtt{torsionMatrix}$ and $\mathtt{magnusMatrix}$.


\bibliographystyle{amsplain}

\end{document}